\documentclass[11pt]{article}
\usepackage{amsmath,amssymb,amsthm,tikz,hyperref}
\usepackage{fancyhdr}

\def\a{\alpha}
\def\b{\beta}

\def\d{\delta}
\def\e{\varepsilon}
\def\z{\zeta}

\def\l{\lambda}

\def\s{\sigma}

\def\w{\omega}

\def\G{\Gamma}

\def\W{\Omega}
\def\AA{\mathbb{A}}

\def\CC{\mathbb{C}}

\def\NN{\mathbb{N}}

\def\PP{\mathbb{P}}

\def\RR{\mathbb{R}}

\def\ZZ{\mathbb{Z}}

\def\cF{\mathcal{F}}
\def\cG{\mathcal{G}}

\def\cO{\mathcal{O}}

\def\cS{\mathcal{S}}

\def\cX{\mathcal{X}}
\def\cY{\mathcal{Y}}

\def\pd{\partial}
\def\ol{\overline}

\def\coker{\mathrm{coker}}

\def\Hom{\mathrm{Hom}}

\def\Op{\cO p}

\newtheorem{theorem}{Theorem}[section]
\newtheorem{corollary}{Corollary}[theorem]
\newtheorem{lemma}[theorem]{Lemma}
\theoremstyle{definition}
\newtheorem{definition}[theorem]{Definition}

\title{Pushforwards of Measures on Real Varieties under Maps with Rational Singularities}
\author{Andrew Reiser}

\begin{document}

\maketitle

\begin{abstract}
Let $X,Y$ be algebraic varieties defined over $\RR$. 
Assume $Y$ is smooth and $X$ is Gorenstein.
Suppose $\varphi:X\to Y$ is a flat $\RR$-morphism such that all the fibers have rational singularities.
We show that the pushforward of any smooth, compactly supported measure on $X$ has a continuous density with respect to any smooth measure with non-vanishing density on $Y$.
This extends a result of Aizenbud and Avni from the $p$-adic case to the archimedean case.
\end{abstract}

\tableofcontents

\section{Introduction}

\begin{definition}
\label{PushforwardDefinition}

Let $X,Y$ be measurable spaces, $\varphi:X\to Y$ a measurable function, and $\mu$ a measure on $X$.
The pushforward of $\mu$ by $\varphi$ is denoted $\varphi_*\mu$ and is defined by $(\varphi_*\mu)(S)=\mu(\varphi^{-1}(S))$ for every $S\subset Y$.

\end{definition}

We are specifically interested in the properties of the pushforward in the case that $X,Y$ are algebraic varieties over $\RR$ equipped with sufficiently nice measures, and $\varphi$ is an $\RR$-morphism of varieties.

\begin{definition}
\label{AbsolutelyContinuous}

Suppose $X$ is a Borel space and $\mu$ is a measure on $X$.
If $\nu$ is another measure on $X$ and for any Borel subset $A\subset X$ we have $\mu(A)=0 \Rightarrow \nu(A)=0$, then $\nu$ is said to be absolutely continuous with respect to $\mu$.

\end{definition}

\begin{definition}
\label{GeneralDensity}

Let $X$ be a measurable space and let $\mu,\nu$ be measures on $X$.
If $\nu$ is absolutely continuous with respect to $\mu$, for any set $A\subset X$ of finite $\mu$-measure, the Radon-Nykodym Theorem says that there is $f_A\in L^1(A,\mu)$ such that $\nu(B)=\int_B f_Ad\mu$ for all measurable subsets $B\subset A$.
Assuming $X$ can be written as a union of sets of finite $\mu$-measure, the functions $f_A$ associated to each set patch together to a globally-defined function $f$ which belongs to $L^1(U)$ for each set $U\in X$ so that $\mu(U)<\infty$, and we call this $f$ the density of $\nu$ with respect to $\mu$.

\end{definition}

We are interested in measures which locally are absolutely continuous with respect to measures coming from algebraic differential forms, where the density function is well-behaved.
The following definition shows how to construct a measure from a differential form.

\begin{definition}
\label{AlgebraicMeasure}

Suppose $X$ is a smooth irreducible	 algebraic variety over $\RR$ and that $\w$ is a rational top differential form on $X$.
We define a measure $|\w|$ on $X(\RR)$ as follows.
Given a relatively-compact open set $U\subset X(\RR)$ and an analytic diffeomorphism $\Psi$ between $U$ and an open subset $W\subset \RR^n$, we may write

$$\Psi^*\w = g dx_1\wedge\cdots\wedge dx_n$$

for some $g:W\to \RR$, and define

$$ |\w|(U) = \int_W |g|d\lambda $$

where $|g|$ is the usual absolute value on $\RR$ and $\lambda$ is the standard Lebesgue measure on $\RR^n$.
By the change of variables formula, this definition is independent of the diffeomorphism $\Psi$.
There is a unique extension of $|\w|$ to a (possibly infinite) Borel measure on $X(\RR)$, which we also denote $|\w|$.

Note that we may assume that $X$ has a smooth real point, as otherwise $X(\RR)=\emptyset$ and there is nothing to do here.

\end{definition}

The following lemma describes some first properties of measures coming from algebraic differential forms.

\begin{lemma}
\label{FirstMeasures}

Suppose $X$ is a smooth, irreducible $\RR$-variety and $\w_1,\w_2$ are two rational top forms on $X$ which are not identically zero.

\begin{itemize}

\item If $\w_1$ is regular, $|\w_1|$ assigns finite values to compact sets $A\subset X(\RR)$.
\item $|\w_1|$ and $|\w_2|$ are absolutely continuous with respect to each other. 
\item If both $\w_1,\w_2$ are regular and $\w_1$ is nowhere-vanishing, the density of $|\w_2|$ with respect to $|\w_1|$ is a continous function on $X(\RR)$.
\item If $\w_1$ has a pole at $x\in X(\RR)$, then $|\w_1|(A)=\infty$ for every open set $A$ containing $x$.

\end{itemize}
\end{lemma}

Some measures on the real points of algebraic varieties which are locally absolutely continuous with respect to measures coming from differential forms are of special interest to us.

\begin{definition}
\label{AlgebraicMeasuresAC}

Suppose $X$ is a smooth, irreducible $\RR$-variety and $\mu$ is a measure on $X(\RR)$.
We say $\mu$ is locally absolutely continuous with respect to algebraic measures at $x$ or is locally (AC) at a point $x\in X(\RR)$ if there exists a Zariski-open neighborhood $U\subset X$ containing $x$ together with a choice of nowhere-vanishing section of the line bundle of top differential forms $\w\in\W_X(U)$ such that when considered as measures on $U(\RR)$, $\mu$ is absolutely continuous with respect to $|\w|$.

\end{definition}

\begin{definition}
\label{AlgebraicMeasures}

Suppose $X$ is a smooth, irreducible $\RR$-variety and $\mu$ is a measure on $X(\RR)$.
We say $\mu$ is locally of continuous density or is locally (CD) at a point $x\in X(\RR)$ if there exists a Zariski-open neighborhood $U\subset X$ containing $x$ together with a choice of nowhere-vanishing section of the line bundle of top differential forms $\w\in \W_X(U)$ such that when considered as measures on $U(\RR)$, $\mu$ is absolutely continuous with respect to $|\w|$ and the density of $\mu$ with respect to $|\w|$ may be chosen to be a continuous function. 

We say $\mu$ has continuous density or is (CD) if $\mu$ is locally (CD) at every point $x\in X(\RR)$.

\end{definition}

Every locally (CD) measure is also locally (AC), and if we have an (AC) measure $\mu$ we wish to show is (CD), it suffices to show that $\mu$ is locally of continuous density at each point. 

\begin{definition}
\label{CompactlySupported}

We call a measure $\mu$ on a measurable space $X$ compactly supported if and only if there exists a compact $K\subset X$ such that for all measurable $A\subset X$, $A\cap (X\setminus K) = \emptyset$ implies that $\mu(A)=0$.

\end{definition}

\begin{definition}
\label{CSCDMeasures}

We call a (CD) measure which is also compactly supported a (CSCD) measure.

\end{definition}

We are interested in finding conditions on a map $\varphi:X\to Y$ between two smooth $\RR$-varieties such that $\varphi_*m$ has continuous density for any (CSCD) measure $m$ on $X(\RR)$.
It is sufficient for $\varphi$ to be a smooth map (see Theorem \ref{WellBehavedSmoothPushforward}), a flat map with all its fibers smooth.
We show that this condition may be relaxed in the following sense: it is enough to let the map be flat with all its fibers having rational singularities (see \ref{RationalSingularitiesDefinition} for a definition of rational singularities).

\begin{definition}
\label{FRS}

Let $X,Y$ be varieties over a field $k$ of characteristic zero, or complex-analytic spaces.
Let $\varphi: X\to Y$ be a morphism.
We say $\varphi$ is (FRS) if $\varphi$ is flat and for all $y\in Y(\ol{k})$, the fiber $X\times_Y y$ is reduced with rational singularities.

\end{definition}

Our main theorem is the following:

\begin{theorem}
\label{IntroTheorem}

Let $X,Y$ be smooth, irreducible varieties over $\RR$ and $\varphi:X\to Y$ be an $\RR$-morphism.
Assume that $\varphi$ is (FRS).
The the pushforward of every (CSCD) measure is again a (CSCD) measure.
\end{theorem}

We prove a stronger version of this theorem:

\begin{theorem}
\label{MainTheorem}

Let $X$ and $Y$ be irreducible varieties over $\RR$ and $\varphi:X\to Y$ be an $\RR$-morphism.
Assume $X$ is Gorenstein and $Y$ is smooth.
Assume that $\varphi$ is (FRS).
The the pushforward of every (CSCD) measure is again a (CSCD) measure.

\end{theorem}

A $p$-adic analogue of this theorem appeared in \cite{AvniAizenbud}.
In proving Theorem \ref{IntroTheorem}, the strategy is similar to the $p$-adic case but with a few detours.
After several reductions, the density function of the pushforward will be shown to be constructible (in this case, log-subanalytic).
From there, the problem will be reduced to the case of a map to a one-dimensional target, and that case will be analyzed.

\subsection{Conventions}

Unless otherwise stated, the following conventions will be obeyed throughout this work:

\begin{itemize}

\item All schemes we consider are finite type over the base field.

\item An algebraic variety is a synonym for a reduced scheme.

\item A morphism of algebraic varieties or schemes means a morphism over the base field.

\item The smooth locus of an algebraic variety $X$ will be denoted by $X^{sm}$.

\item The smooth locus of a morphism of schemes $\varphi: X\to Y$ will be denoted by $X^S$.

\item When referring to a point, we mean a closed point over the base field.

\item For a variety $X$ defined over $\RR$ or $\CC$, we will often consider the set of $\RR$ or $\CC$ points together with the analytic topology, and we use $X(\RR)$ or $X(\CC)$ to denote these topological spaces.

\item When referring to geometric spaces, we will use Roman letters for schemes and script letters for analytic spaces.

\item For an algebraic or analytic variety, we denote the dualizing complex $\W^\bullet$, and in the case that this complex is concentrated in a single degree, we use $\W$ for the dualizing sheaf.
We identify the restriction of the dualizing sheaf to the smooth part of the variety with top differential forms on the smooth locus.

\item If $A=(a_1,\cdots,a_n)$ is a sequence of integers or real numbers, we write $A>c$ (respectively, $A\geq c$) if each $a_i>c$ (respectively $a_i\geq c$).

\end{itemize}

\section{General Properties of the Pushforward}

Let $\varphi:X\to Y$ be a smooth map of smooth $k$-varieties. 
Denote by $\W_{X/Y}$ the sheaf of relative top differential forms on $X$ with respect to $Y$.
Since $\varphi$ is smooth, $\W_{X/Y}$ is an invertible sheaf on $X$ and there is an isomorphism $\W_X\to \varphi^*\W_Y\otimes \W_{X/Y}$ such that for every field extension $K\supset k$ and every point $x\in X(K)$, the isomorphism of fibers $$\wedge^{\dim X}T_x^*X \to \left(\wedge^{\dim Y} T_{\varphi(x)}^*Y\right)\otimes\left(\wedge^{\dim \varphi^{-1}(\varphi(x))}T_{x}^*\varphi^{-1}(\varphi(x))\right)$$ comes from the short exact sequence of vector spaces $$ 0\to T_x\varphi^{-1}(\varphi(x))\to T_xX\stackrel{d\varphi}{\to} T_{\varphi(x)}Y\to 0.$$

For a top form $\w_X\in \G(X,\W_X)$ and a nowhere-vanishing top form $\w_Y\in \G(Y,\W_Y)$, there exists a unique element $\eta\in \G(X,\W_{X/Y})$ such the the image of $\eta\otimes \varphi^*\w_Y$ under the isomorphism $\W_{X/Y}\otimes \varphi^*\W_Y\to \W_X$ is $\w_X$. 

\begin{definition}
\label{QuotientDifferentialForms}

For the rest of the paper, we will refer to the $\eta$ defined above as $\frac{\w_X}{\varphi^*\w_Y}$. 
In cases where $\varphi: X\to Y$ is not smooth, we shall carry out the construction on the smooth locus of $\varphi$.

\end{definition}


Next, we discuss how to extend the concept of a locally (CD) measure (definition \ref{AlgebraicMeasure}) to a variety with mild singularities.

\begin{definition}
\label{AlgebraicMeasuresGeneralization}

Suppose $X$ is an irreducible Gorenstein algebraic variety defined over $\RR$ with dualizing line bundle $\W_X$.
For $\w\in\W_X(X)$, we define a measure $|\w|$ on $X(\RR)$ as follows.
Recall that $\W_X|_{X^{sm}}$ is isomorphic to the line bundle of top differential forms on $X^{sm}$.
Given a relatively-compact open set $U\subset X(\RR)$ and an analytic diffeomorphism $\Psi$ between $U\cap X^{sm}(\RR)$ and an open subset $W\subset \RR^n$, we may write $$\Psi^*\w = gdx_1\wedge \cdots\wedge dx_n$$ for some $g:W\to\RR$, and define $$|\w|(U)=\int_W |g|d\l$$ where $|g|$ is the usual absolute value on $\RR$ and $\l$ is the standard Lebesgue measure on $\RR^n$.
By the change of variables formula, this definition is independent of the diffeomorphism $\Psi$.
There is a unique extension of $|w|$ to a (possibly infinite) Borel measure on $X(\RR)$, which we also denote $|\w|$.


\end{definition}

With this modification, each of Definitions \ref{AlgebraicMeasuresAC}, \ref{AlgebraicMeasures}, \ref{CompactlySupported} and \ref{CSCDMeasures} generalize in a natural way to the case of $X$ Gorenstein, so we may speak of (AC), (CD), and (CSCD) measures on Gorenstein varieties.


\begin{theorem}
\label{WellBehavedSmoothPushforward}

Let $\varphi:X\to Y$ be a smooth map of smooth varieties over $\RR$.

\begin{enumerate}

\item The pushfoward of a compactly supported measure with continuous density is compactly supported with continuous density.

\item Let $\w_X,\w_Y$ be top differential forms on $X,Y$ respectively. 
Assume $\w_Y$ is nowhere vanishing and $f$ is a compactly supported continuous function on $X(\RR)$. 
The measure $\varphi_*(f|\w_X|)$ is absolutely continuous with respect to $|\w_Y|$ and its density at a point $y\in Y(\RR)$ is equal to $\int_{\varphi^{-1}(y)(\RR)} f\left| \frac{\w_X}{\varphi^*\w_Y}|_{\varphi^{-1}(y)} \right|$

\end{enumerate}

\end{theorem}

\begin{proof}

First, we prove that the pushforward of a compactly supported measure is compactly supported. 
Suppose $\mu$ is a compactly-supported measure on $X(\RR)$.
If $K\subset X(\RR)$ is a compact set so that for all measurable $A\subset X(\RR)\setminus K$, we have $\mu(A)=0$, then $\varphi(K)\subset Y(\RR)$ is a compact set so that for all measurable $B\subset Y(\RR)\setminus\varphi(K)$, we have $(\varphi_*\mu)(B)=0$.

To prove that the pushforwards of a (CSCD) measure is again (CSCD), we prove the formula for the density first.
In order to prove the second claim, we may use a partition of unity and an analytic change of coordinates to reduce to the case where $X$ and $Y$ are open subsets of $\AA_\RR^n$ and $\AA_\RR^m$, respectively, and the morphism is linear projection from $\AA_\RR^n\to\AA_\RR^m$, where $\AA_\RR^m\subset\AA_\RR^n$ is embedded as the subspace of the final $m$ coordinates.
Let $A\subset Y(\RR)$ be a measurable subset.
Let $1_{S}$ be the characteristic function of a subset $S$.
We may fix functions $\Psi,\Xi$ with $\Xi$ is nonvanishing on $Y$ so that $\w_X=\Psi dx_1\cdots dx_n$, $\w_Y=\Xi dx_{n-m+1}\cdots dx_{n}$ and $\frac{\w_X}{\varphi^*\w_Y}=\frac{\Psi}{\Xi} dx_1\cdots dx_{n-m}$.
Rewrite $$\varphi_*(f|\w_X|)(A)=(f|\w_X|)(\varphi^{-1}(A))=\int_{X(\RR)} 1_{\varphi^{-1}(A)}f|\w_X|=\int_{\RR^n} 1_{X(\RR)\cap \RR^{n-m}\times(A)} f|\Psi| dx_1\cdots dx_n .$$

By an application of Fubini's Theorem, this may be rewritten as $$\int_{\RR^m}\left(\int_{\RR^{n-m}} 1_{X\cap (A\times \RR^{n-m})} f\left|\frac{\Psi}{\Xi}\right| dx_1\cdots dx_{n-m}\right) |\Xi|dx_{n-m+1}\cdots dx_{n},$$ which exhibits $\varphi_*\mu$ as absolutely continuous with respect to $|\w_Y|$ and also gives the formula for the density.
This formula, along with the observation that $f|\frac{\Psi}{\Xi}|$ is uniformly continuous, immediately implies that $\varphi_*\mu$ is (CD).

\end{proof}

\begin{corollary}
\label{LocallyDominantPushforward}


Let $\varphi:X\to Y$ be a locally dominant map between smooth $\RR$-varieties and denote the smooth locus of $\varphi$ by $X^S$.

\begin{enumerate}
\item If $m_X$ is a (CSCD) measure on $X(\RR)$, then $\varphi_*m_X$ is a compactly supported measure on $Y$ which is (AC).
\item In particular, if $m_X=f|\w_X|$ for a top differential form $\w_X\in\W_{X}(X)$ and continuous, compactly supported function $f:X(\RR)\to \RR$ and $\w_Y$ is a nowhere vanishing top differential form on $Y$, $\varphi_*m_X$ is absolutely continuous with respect to $|\w_Y|$ and has density given by 

$$y\mapsto \int_{(X^S\cap\varphi^{-1}(y))(\RR)} f \left| \frac{\w_X}{\varphi^*\w_Y}|_{X^S\cap\varphi^{-1}(y)} \right|$$

as a function from $Y(\RR)$ to $\RR_{\geq 0}\cup \{\infty\}$.
\end{enumerate}

\end{corollary}

\begin{proof}

Since the claims above are local on $Y$, we may assume that $Y$ is affine.
We may cover $X$ by affine opens $U_i$ and construct a partition of unity on $X(\RR)$ subordinate to the cover $U_i(\RR)$.
With this partition of unity, we may assume that $X$ is affine and choose an embedding of it in to $\AA^N$ for some $N$.
Let $X(\RR)^{s,n}$ be the set of points of $X^S(\RR)$ such that their distance in the natural metric on $\AA^N(\RR)$ to the set $(X\setminus X^S)(\RR)$ is at least $\frac{1}{n}$, and let $g_n$ be a $C^\infty$ function supported on $X(\RR)^{s,n+1}$ which takes the value $1$ on $X(\RR)^{s,n}$ and satisfies $0\leq g_n\leq 1$ everywhere.

Since $g_nm_X$ satisfies the assumptions of Theorem \ref{WellBehavedSmoothPushforward} on $X^S(\RR)$, we may apply Theorem \ref{WellBehavedSmoothPushforward} to see that $\varphi_*(g_nm_X)$ is absolutely continuous with respect to $m_Y$, so it has a density, denoted $h_n$. 
Since $h_n$ is nondecreasing in $n$ and $\int_Y h_nm_Y \leq m_X(X(\RR))<\infty$, we have by Lebesgue's Monotone Convergence Theorem that $h:=\lim_{n\to\infty} h_n$ exists and belongs to $L_1(Y(\RR),m_Y)$. 
Since $\varphi$ is locally dominant and $m_X$ is (CD), we have that $\varphi_*(m_X)=(\varphi|_{X^S})_*(m_X|_{X^S})$. 
After integrating against continuous functions and another application of the Monotone Convergence Theorem, we see that $\varphi_*(m_X)=(\varphi|_{X^S})_*(m_X|_{X^S})=hm_Y$, and the latter is absolutely continuous with respect to $m_Y$.

The statement on the integral formula representation for the density is clear from Theorem \ref{WellBehavedSmoothPushforward}.

\end{proof}

Next, there are some lemmas about how (CD) measures and their pushforward interact with rational singularities. 
For a reminder on the definition and equivalent characterizations of rational singularities, see Section A.3 in Appendix A, especially Definition \ref{RationalSingularitiesDefinition} and Theorems \ref{MoreRationalSingularities} and \ref{MoreRationalSingularitiesAnalytic}. 

\begin{lemma}
\label{RationalSingFinIntegral}

Let $V\subset \AA^n$ be an $\RR$-variety with rational singularities. 
Let $U\subset V$ be an open smooth subset such that $V\setminus U$ has codimension two in $V$. 
Additionally, let $\w$ be a regular top differential form on $U$. 
Then for any compact subset $K\subset \AA^n(\RR)$, the integral $\int_{U(\RR)\cap K} |\w|$ is finite. 

\end{lemma}

\begin{proof}



Take a strong resolution of singularities $\pi:\widetilde{V}\to V$. 
By Theorem \ref{MoreRationalSingularities} part (7), there exists a top differential form $\widetilde{\w}$ on $\widetilde{V}$ so that the restriction of $\widetilde{\w}$ to $X^{sm}$ is equal to $\w$.
In particular, $\widetilde{\w}|_U=\w$.
Computing, we see

$$ \int_{U(\RR)\cap K} |\w| = \int_{\pi^{-1}(U(\RR)\cap K)} |\widetilde{\w}| = \int_{\pi^{-1}(V(\RR)\cap K)} |\widetilde{\w}|$$

where the final integral is a continuous function over a compact set and is therefore finite.

\end{proof}

This shows that an algebraic measure on the real points of a Gorenstein variety with rational singularities assigns finite values to compact sets.

\begin{lemma}
\label{AbsolutelyContinuousUpgraded}


Let $X$ be a Gorenstein variety with rational singularities over $\RR$. 
Let $Y$ be a smooth variety over $\RR$. 
Let $\varphi:X\to Y$ be a locally dominant map, and let $X^{s}$ be the smooth locus of $\varphi$. 

\begin{enumerate}
\item If $m_X$ is a (CSCD) measure on $X(\RR)$, then $\varphi_*m_X$ is a compactly supported measure on $Y$ which is (AC).
\item In particular, if $m_X=f|\w_X|$ for a section $\w_X\in\W_{X}(X)$ and continuous, compactly supported function $f:X(\RR)\to \RR$ and $\w_Y$ is a nowhere vanishing top differential form on $Y$, $\varphi_*m_X$ is absolutely continuous with respect to $|\w_Y|$ and has density given by 

$$y\mapsto \int_{(X^S\cap\varphi^{-1}(y))(\RR)} f \left| \frac{\w_X}{\varphi^*\w_Y}|_{X^S\cap\varphi^{-1}(y)} \right|$$

as a function from $Y(\RR)$ to $\RR_{\geq 0}\cup \{\infty\}$.
\end{enumerate}




\end{lemma}

\begin{proof}

As in \ref{LocallyDominantPushforward}, it suffices to prove the second claim assuming that $X,Y$ are affine.
Choose a resolution of singularities $\pi: \widetilde{X}\to X$, let $\widetilde{\varphi}=\varphi\circ\pi$, and let $\w_{\widetilde{X}}$ be a top differential form that coincides with $\w_X$ on an open dense set. 
By \ref{LocallyDominantPushforward}, we see that the measure $\widetilde{\varphi}_*((f\circ\pi)|\w_{\widetilde{X}}|)=\varphi_*(f|\w_X|)$ is absolutely continuous with respect to $|\w_Y|$ and has density given by

$$y\mapsto \int_{(X^S\cap\varphi^{-1}(y))(\RR)} f \left| \frac{\w_X}{\varphi^*\w_Y}|_{X^S\cap\varphi^{-1}(y)} \right|$$

Let $U\subset X$ be an open dense set such that $\pi|_U$ is an isomorphism. 
Then $\pi^{-1}(U)$ is open dense in $\widetilde{X}$, and there exists an open dense subset $V\subset Y$ such that for any $y\in V$, $U\cap\varphi^{-1}(y)$ is dense in $\varphi^{-1}(y)$ and the set $\pi^{-1}(U)\cap\widetilde{\varphi}^{-1}(y)$ is dense in $\widetilde{\varphi}^{-1}(y)$. 
Passing from $(X^{s}\cap\varphi^{-1}(y))(\RR)$ as our domain of integration to $(U\cap X^{s}\cap\varphi^{-1}(y))(\RR)$, then pulling back to $\widetilde{X}$ and using the fact that $\pi^{-1}(U)$ is open dense, we obtain the result.

\end{proof}

\section{O-Minimal Geometry}

In this section, we establish all of the results in o-minimal geometry we will need for the rest of the paper.
Subsection one gives the basic definitions. 
Subsection two discusses the continuity of functions definable in an o-minimal structure and establishes that continuity of log-subanalytic functions may be checked along analytic curves.
Subsection three deals with approximating arbitrary functions and sets by semialgebraic functions and sets.
Subsection four deals with the integration of definable functions.

\subsection{First Definitions}

\begin{definition}

An o-minimal structure on the field $(\RR,+,\cdot,0,1)$ is a collection $\cS=\{\cS_n\}_{n\in \ZZ_{\geq 0}}$, where each $\cS_n$ is a family of subsets of $\RR^n$ satisfying the following properties:

\begin{enumerate}
\item $\cS_n$ contains all algebraic subsets of $\RR^n$;
\item $\cS_n$ is a Boolean subalgebra of the power set of $\RR^n$, i.e. it is closed under finite intersection, finite union, and complement;
\item If $A\in\cS_n$, $B\in\cS_m$, then $A\times B\in \cS_{m+n}$;
\item If $\pi: \RR^n\times\RR \to \RR^n$ is the natural projection and $A\in \cS_{n+1}$, then $A\in \cS_n$;
\item $\cS_1$ is exactly the collection of finite unions of points and intervals.
\end{enumerate}

\end{definition}

\begin{definition}
\label{SemialgebraicDefinition}

A subset $A\subset \RR^n$ is semialgebraic if it may be described as $$A = \bigcup_{i=1}^{p} \bigcap_{j=1}^q \{x\in \RR^n \mid f_{ij}(x) \ast 0\}$$

where $f_{ij}(x)$ are polynomials and $\ast$ stands for any of the symbols $>,<,=$. 

Let $A\subset \RR^n$ be semialgebraic. A function $f:A\to \RR^m$ is semialgebraic if its graph is a semialgebraic subset of $\RR^{n+m}$.

\end{definition}

\begin{theorem}
\label{SemialgebraicProperties}

Semialgebraic sets form an o-minimal structure, and every o-minimal structure contains the semialgebraic sets.
In addition to satisfying the conditions of an o-minimal structure, semialgebraic sets and functions possess the following properties:

\begin{enumerate}
\item The semialgebraic sets are exactly the definable subsets of $\RR^n$ in the structure $\mathcal{R}=(\RR,0,1,+,\cdot)$ (equivalently, $\mathcal{R}=(\RR,0,1,+,\cdot,\leq)$);
\item Each semialgebraic set has a finite number of connected components, each of which is semialgebraic;
\item If $A$ is semialgebraic, then the closure, interior, and boundary are also semialgebraic;
\item The sum, product, and scalar multiples of semialgebraic functions are semialgebraic;
\item The Euclidean distance function from a point to a fixed non-empty semialgebraic set is semialgebraic.
\end{enumerate}

\end{theorem}

\begin{proof}

See, for example, section 1.1 of \cite{Denkowska}.

\end{proof}

We will require greater generality than semialgebraic sets give us access to. In order to access the definition of globally subanalytic sets, we require two intermediate definitions.

\begin{definition}

Let $M$ be a real analytic manifold.
A subset $A\subset M$ is called semianalytic if for any $x\in M$ there exists an open neighborhood $U\ni x$ and analytic functions $f_i$, $g_{ij}$ on $U$ such that 

$$A\cap U = \bigcup_{i=1}^p\bigcap_{j=1}^q \{x\in U\mid f_{i}(x)=0, g_{ij}(x)>0\}$$

\end{definition}

\begin{definition}

Let $M$ be a real-analytic manifold.
A subset $E\subset M$ is called subanalytic if for any $x\in M$ there exists an open neighborhood $U\ni x$ together with a real analytic manifold $N$ and a semianalytic set $A\subset M\times N$ such that $E\cap U =\pi(A)$, where $\pi:M\times N\to M$ is the natural projection.

\end{definition}

\begin{definition}

A subset $E\subset \RR^n$ is called globally subanalytic if its image under the natural embedding of $\RR^n\hookrightarrow \PP^n$ given by $(x_1,\cdots,x_n)\mapsto (1:x_1:\cdots:x_n)$ is subanalytic. 

Let $E\subset \RR^n$ be globally subanalytic.
A function $f:E\to \RR^m$ is globally subanalytic if its graph is a globally subanalytic subset of $\RR^{n+m}$.

\end{definition}

\begin{theorem}
\label{GloballySubanalyticProperties}

Globally subanalytic sets and functions have the following properties:

\begin{enumerate}

\item Globally subanalytic sets are exactly the sets definable in the expansion of the real field by all restricted analytic functions $f:\RR^n\to \RR$, i.e. functions of the form $1_{[-1,1]^n}f_0$, where $f_0$ is an analytic function defined on an open neighborhood of $[-1,1]^n$;
\item Globally subanalytic sets form an o-minimal structure;
\item The interior, closure, and boundary of a globally subanalytic set are globally subanalytic;
\item The connected components of a globally subanalytic set are each globally subanalytic and the collection of connected components is finite.

\end{enumerate}

\end{theorem}

\begin{proof}

See, for example, section 2.3 of \cite{Denkowska}.

\end{proof}

\begin{definition}
\label{LogSubanalytic}

A log-subanalytic function on a globally subanalytic set $E\subset\RR^n$ is a function which may be written in either of the two following equivalent forms:

\begin{enumerate}
\item as a polynomial in a finite number of globally subanalytic functions and logarithms of strictly positive globally subanalytic functions;
\item as a polynomial in a finite number of strictly positive globally subanalytic functions and their logarithms.
\end{enumerate}

Note that the usual logarithm is excluded from this class of functions.

\end{definition}

\subsection{Continuity}


The following theorem shows that functions definable in a given o-minimal structure are piecewise continuous, and establishes a convenient description of the pieces.

\begin{theorem}
\label{DefinablyPiecewiseContinuous}

Let $X\subset \RR^n$ and $Y\subset \RR^m$ be definable sets in a fixed o-minimal structure $\cS$.
Let $f:X\to Y$ be a function definable in $\cS$.
There exists a finite definable partition of $X=X_1\cup\cdots\cup X_n$ such that $f|_{X_i}$ is continuous for each $i$.

\end{theorem}

\begin{proof}

See Theorem 2.12 of \cite{CosteOMinimal}.

\end{proof}

Subanalytic subsets of manifolds may also be treated with the tools of resolution of singularities, as in Hironaka's Rectilinearization Theorem \cite{BierstoneMillman}:

\begin{theorem}
\label{HironakaRectilinearization}

Assume that $M$ is a pure-dimensional real analytic manifold and $E\subset M$ a subanalytic subset. 
Let $K$ be a compact subset of $M$. 
Then there are finitely many real analytic mappings $\psi_i:\RR^m\to M$ such that:

\begin{enumerate}
\item For each $i$, there is a compact subset $L_i$ of $\RR^m$, such that $\bigcup_i \psi_i(L_i)$ is a neighborhood of $K$ in $M$;
\item For each $i$, $\psi_i^{-1}(E)$ is a union of quadrants in $\RR^m$.
\end{enumerate}

Where a quadrant $Q\subset\RR^m$ is a set described by the $m$ equations (one for each $1\leq i\leq m$) $x_i\ast 0$, where $\ast$ is any of the symbols $>,=,<$.

\end{theorem}

The combination of rectilinearization and piecewise continuity give the following theorem, which says that a log-subanalytic function from a manifold to $\RR$ is continuous if and only its restriction to every analytic curve is continuous.

\begin{theorem}
\label{LogSubanalyticDiscontinuity}

Let $M\subset \RR^n$ be a real-analytic manifold of dimension $m$ and $f:M\to \RR$ be a log-subanalytic function.
If $f$ is discontinuous at a point $x_0\in M$, then there exists an analytic map $g:\RR\to M$ such that $x_0\in g(\RR)$ and $f|_{g(\RR)}$ is discontinuous.

\end{theorem}

\begin{proof}

Since $f:M\to\RR$ is log-subanalytic, we may pick strictly positive globally subanalytic functions $\varphi_1,\cdots,\varphi_r:M\to\RR_{>0}$ and a polynomial $p\in\RR[x_1,\cdots,x_r,y_1,\cdots,y_r]$ such that $f=p(\varphi_1,\cdots,\varphi_r,\log\varphi_1,\cdots,\log\varphi_r)$.
By Theorem \ref{DefinablyPiecewiseContinuous}, for each $\varphi_i$ we may decompose $M$ into $n_i$ globally subanalytic pieces $M_{i,1},\cdots,M_{i,n_i}$ such that $\varphi_i$ is continuous on each $M_{i,j}$.

For each $\alpha=(\alpha_1,\cdots,\alpha_r)\in \prod_{1\leq i\leq r} \{1,2,\cdots,n_i\}=A$, let $M_{\alpha}=\bigcap_{1\leq i\leq r} M_{i,\alpha_i}$.
Then each $M_\alpha$ is globally subanalytic, the $M_{\alpha}$ partition $M$, and the restriction of $f$ to each $M_{\alpha}$ is continuous.
One of the following three claims must hold:

\begin{enumerate}
\item $f$ extends to a continuous function on $\ol{M_{\alpha}}$ for all $\alpha\in A$;
\item There exists $\alpha_0\in A$ and a point $m_0\in\ol{M_{\alpha_0}}$ such that $\lim_{m\in M_{\alpha_0},m\to m_0} f(m)=\pm\infty$;
\item There exist $\alpha_1,\alpha_2\in A$ such that $\ol{M_{\alpha_1}}\cap\ol{M_{\alpha_2}}\neq \emptyset$ and $\inf_{m_1\in M_{\alpha_1},m_2\in M_{\alpha_2}} |f(m_1)-f(m_2)|\neq 0$.
\end{enumerate}

If the first claim holds, then $f:M\to\RR$ is continuous and we have nothing to prove.

If the second claim holds, note that $M_{\alpha}$ is positive-dimensional, as otherwise it is finite and the claim cannot hold.
We may apply Theorem \ref{HironakaRectilinearization} to $M_{\alpha_0}\subset M$ to get an analytic map $h:\RR^{\dim{M}}\to M$ such that $h^{-1}(M_{\alpha_0})$ is a union of quadrants, at least one of which is positive dimensional.
Since $m_0\notin M_{\alpha_0}$, $h^{-1}(m_0)\notin h^{-1}(M_{\alpha_0})$.
On the other hand, $m_0\in\ol{M_{\alpha_0}}$ so $h^{-1}(m_0)\in\ol{h^{-1}(M_{\alpha_0})}$.
This means we can pick a line $L_0\subset \RR^{\dim M}$ such that $h^{-1}(M_{\alpha_0})\cap L_0$ is an open ray which contains a point of $h^{-1}(m_0)$ in it's closure.
After identifying $L_0$ and $\RR$, the map $h|_{L_0}:L_0\to M$ satisfies the conclusions of the lemma.

If the third claim holds, pick a point $m_0\in\ol{M_{\alpha_1}}\cap\ol{M_{\alpha_2}}$.
We may assume both $M_{\alpha_1}$ and $M_{\alpha_2}$ are positive-dimensional, as otherwise one of them is finite and we may apply the same logic from the previous paragraph.
After possibly subdividing the partition, we may assume $m_0\notin M_{\alpha_1}$ and $m_0\notin M_{\alpha_2}$. 
Apply Theorem \ref{HironakaRectilinearization} to $M_{\alpha_1},M_{\alpha_2}\subset M$ to get analytic maps $h_1,h_2:\RR^{\dim{M}}\to M$ such that each $h_i^{-1}(M_{\alpha_i})$ is a union of quadrants, and at least one of these quadrants is positive dimensional for each $i$.
By our assumption on $m_0$, we have that $h_i^{-1}(m_0)\cap h_i^{-1}(M_{\alpha_i})=\emptyset$ for $i=1,2$.
On the other hand, $m_0\in\ol{M_{\alpha_i}}$, so $h^{-1}(m_0)\subset\ol{h_i^{-1}(M_{\alpha_i})}$ for both $i$.
This means we can pick a lines $L_1,L_2\subset \RR^{\dim M}$ such that $h_i^{-1}(M_{\alpha_I})\cap L_i$ is an open ray which contains a point of $h^{-1}(m_0)$ in it's closure.
Then $h_1(L_1)$ and $h_2(L_2)$ are two analytic curves in $M$ which meet at $m_0$.
If $f$ was continuous when restricted to both $h_1(L_1)$ and $h_2(L_2)$, then it would be continuous on their union, which contradicts our assumption that $\inf_{m_1\in M_{\alpha_1},m_2\in M_{\alpha_2}} |f(m_1)-f(m_2)|\neq 0$.
So, without loss of generality, $f$ is discontinuous when restricted to $h_1(L_1)$.
After identifying $L_1$ and $\RR$, the map $h_1|_{L_1}:L_1\to M$ satisfies the conclusions of the lemma.

\end{proof}

\subsection{Approximation by Semialgebraic Functions}

Firstly, we recall the Stone-Weierstrass Theorem:

\begin{theorem}
\label{StoneWeierstrass}

Let $K$ be a compact Hausdorff topological space.
Let $C(K,\RR)$ be the space of real-valued continuous functions on $S$.
Let $A\subset C(K,\RR)$ be a subalgebra.
$A$ is dense in $C(K,\RR)$ in the topology of uniform convergence if and only if $A$ contains a nonzero constant function and $A$ separates points.

\end{theorem}

We will show that the semialgebraic continuous functions on any compact semialgebraic set fulfill the conditions of the theorem.

\begin{lemma}
\label{SemialgebraicDensity}

Fix $K$ a semialgebraic compact subset of $\RR^n$.
Let $A$ be the collection of continuous semialgebraic functions $K\to \RR$.
The following are true:

\begin{enumerate}
\item $A$ is a subalgebra of $C(K,\RR)$.
\item $A$ separates points.
\item $A$ contains a nonzero constant function.
\end{enumerate}

\end{lemma}

\begin{proof}

From Theorem \ref{SemialgebraicProperties} on the properties of semialgebraic functions, they form an algebra.
To show that continuous semialgebraic functions separate points, recall that the function which returns the Euclidean distance to a point is semialgebraic, so $d(x_1,x)$ is a continuous semialgebraic function which separates any two distinct points $x_1,x_2$.
Since $K$ is semialgebraic, every constant function on it is semialgebraic, and in particular, the constant function $1_K$ is in $A$.

\end{proof}

Next, we will show that for any affine variety $X(\RR)$ and any compact set $K\subset X(\RR)\subset\RR^n$, there exists a continuous compactly supported semialgebraic cutoff function $b_K:X(\RR)\to \RR$ taking the value $1$ on $K$.

\begin{lemma}
\label{SemialgebraicCompactTrick}

Suppose $X$ is an affine variety defined over $\RR$. 
Then the embedding $X\hookrightarrow \AA^n_\RR$ gives an inclusion $X(\RR)\subset \RR^n$.
Suppose $K\subset X(\RR)$ is a compact set.
Then there exists a continuous, compactly supported, semialgebraic function $b_K:X(\RR)\to \RR$ such that $0\leq b_K\leq 1$ and $b_K|_K=1$.

\end{lemma}

\begin{proof}

As $X(\RR)$ is locally compact, for every point $x\in X(\RR)$ there exists a compact neighborhood $N_x\subset X(\RR)$ containing $B(x,q_x)\cap X(\RR)$ for some real number $q_x>0$.
For each $x\in X(\RR)$, let $\psi_x(y) = \min(0,\frac{q_x}{2}-d(x,y))$ where $d$ is the induced metric from the inclusion $X(\RR)\hookrightarrow \RR^n$.
Then each $\psi_x:X(\RR)\to \RR$ is continuous, semialgebraic, and compactly supported (its support is a closed subset of a compact subset, $N_x$).
The sets $\psi_x^{-1}((0,\infty))$ cover $X(\RR)$, and since $\operatorname{Supp}(f)$ is compact, there exists a finite set $\{x_1,\cdots,x_N\}\subset X(\RR)$ such that $\cup_{1\leq i\leq N} \psi_{x_i}^{-1}((0,\infty))$ cover $\operatorname{Supp}(f)$.
As a result, $\sum_{i=1}^N \psi_{x_i}$ is a compactly supported continuous semialgebraic function taking positive values on $\operatorname{Supp}(k)$.
After multiplication by some large enough $C>0$, $b_K:=\min(1,C\sum_{i=1}^N \psi_{x_i})$ is a semialgebraic, continuous, compactly supported function taking the value $1$ on $\operatorname{Supp}(f)$.

\end{proof}

This implies that for any $K$, $b_K^{-1}(1)$ is a compact semialgebraic subset containing $K$.
In particular, any compactly supported function $f:X(\RR)\to\RR$ has a compact semialgebraic subset of $X(\RR)$ which contains its support.
This observation, combined with Theorem \ref{StoneWeierstrass} and Lemma \ref{SemialgebraicDensity} gives the following corollary.

\begin{corollary}
\label{SWPunchline}

Suppose $X\subset \AA^n_\RR$ is an affine $\RR$-variety.
Let $f:X(\RR)\to\RR$ be a nonnegative compactly-supported continuous function.
Then there exists a sequence of nonnegative continuous compactly supported semialgebraic functions $f_n:X(\RR)\to\RR$ which converges uniformly to $f$ and at all points, $f_n\leq f$.

\end{corollary}

\begin{proof}

By Lemma \ref{SemialgebraicCompactTrick}, we may fix a compact semialgebraic subset $K\subset \RR^n$ containing $\operatorname{Supp} f$.
By an application of Theorem \ref{StoneWeierstrass} and Lemma \ref{SemialgebraicDensity}, we may find a continuous semialgebraic function $g_{n}$ such that $||(f-\frac{1}{2n})-g_n|| < \frac{1}{2n}$ on $K$.
For each $n$, $\{x\in K \mid g_n(x)> 0\}$ is a semialgebraic subset of $K$ and therefore of $X(\RR)$.
Define $f_n=\max(g_n,0)$, which is semialgebraic.
Furthermore, $f_n\geq 0$ at all points, and $f_n$ is continuous.
Since $f(x)=f_n(x)=0$ for all $x\in X(\RR)\setminus K$ and $|f(x)-f_n(x)|\leq |f(x)-g_n(x)|\leq \frac{1}{n}$ for all $x\in K$ and all $n\geq 1$, we see that $f_n\to f$ uniformly as $n\to \infty$.

\end{proof}

\subsection{Integration}

In this section, we discuss the parameterized integration of semialgebraic functions on semialgebraic domains.

\begin{definition}
\label{InfiniteIntegral}

Suppose $X\subset \RR^n$ is a semialgebraic set and let $f:X\times \RR^m\to \RR$ be a semialgebraic function.
Let $\infty(X,f) = \{x\in X \mid \int_{\RR^m} |f(x,\xi)| d\xi =\infty\}$ be the locus in $X$ where $f$ is not integrable over the fibers.

\end{definition}

\begin{lemma}
\label{InfiniteLocusExcision}

Suppose $X\subset \RR^n$ is a semialgebraic set and let $f:X\times \RR^m\to \RR$ be a semialgebraic function.
Then $f$ is measurable and $\infty(X,f)\subset X$ is semialgebraic.

\end{lemma}

\begin{proof}

See Kaiser\cite{KaiserTame} Proposition 1.1 and Theorem 2.3a.

\end{proof}

The follwing theorem of Cluckers and Miller \cite{CluckersMiller} implies that the integral of any log-subanalytic function is again log-subanalytic.

\begin{theorem}
\label{CluckersMillerTheorem}

Let $E\subset\RR^n$ be a globally subanalytic set, $m\geq 0$ an integer, and $f:E\times \RR^m\to \RR$ a log-subanalytic function such that $f(x,-)$ is Lebesgue integrable for all $x\in X$. 
Then $I_X(f):X\to \RR$ given by $x\mapsto \int_{\RR^m} f(x,y)dy$ is log-subanalytic.

\end{theorem}

If $f$ and $X$ from Theorem \ref{CluckersMillerTheorem} are actually semialgebraic, one may semialgebraically alter $f$ to remove the locus where it is not integrable by an application of Lemma \ref{InfiniteLocusExcision}.
In this way, one may make sense of $I_X(f)$ for arbitrary semialgebraic $f$ by declaring $I_X(f)$ to return $\infty$ on $\infty(X,f)$.

In order to apply Theorem \ref{CluckersMillerTheorem}, we will establish some results similar to the techniques of appendix A in \cite{AvniAizenbud}.

\begin{lemma}
\label{SemialgebraicFiniteFibers}

Suppose that $X\subset \AA^N_\RR$ and $Y\subset \AA^M_\RR$ are affine $\RR$-varieties, that $\pi:X\to Y$ is a morphism with finite fibers, and that $f:X(\RR)\to\RR$ is semialgebraic.
Then the function $y\in Y(\RR)\mapsto \sum_{x\in\pi^{-1}(y)(\RR)} |f(x)|$ is semialgebraic.
\end{lemma}

\begin{proof} 
We will show that there are semialgebraic functions $g_i:Y(\RR)\to \RR$ such that for every $y\in Y(\RR)$ we have $\{|f(x)|: x\in\pi^{-1}(y)(\RR) \wedge f(x)\neq 0\} \subset \{g_1(y),\cdots,g_m(y)\}$ and, after passing to a semialgebraic partition, the number of points $x\in\pi^{-1}(y)(\RR)$ such that $|f(x)|=g_i(y)$ is a constant $n_i$. 
This implies that $\sum_{x\in\pi^{-1}(y)(\RR)} |f(x)|=\sum n_ig_i(y)$ on each piece, and therefore that both sums are semialgebraic functions.

Let $N$ be the maximum size of the fiber of $\pi$. 
We will show by induction on $i$ that there exists a finite semialgebraic partition of $Y(\RR)$ and on each piece $A$ there are semialgebraic functions $g_1,\cdots,g_i$ and natural numbers $n_1,\cdots,n_i$ such that for each $y\in A$ the following hold:

\begin{enumerate}
\item $\#\{x\in \pi^{-1}(y): |f(x)|=g_i(y)\} = n_i$
\item $\#\{x\in \pi^{-1}(y): f(x)\neq 0 \wedge |f(x)|\notin \{g_1(y),\cdots,g_i(y)\}\} \leq N-i$
\end{enumerate}

If we can show this claim for $i=N$, then we will have completed the proof.

Suppose that $g_j$ and $n_j$ were chosen for all $j<i$, and let 

$$g_i(y) = \min \{|f(x)| : \pi(x)=y \wedge |f(x)|\neq g_1(y),\cdots,g_{i-1}(y)\}$$

As the minimum of a finite number of semialgebraic functions is semialgebraic, $g_i$ is semialgebraic. 
For each $n$, the set of $y\in Y$ such that there are at least $n$ distinct points $x_k$ in the fiber over $y$ such that $|f(x_k)|=g_i(y)$ is semialgebraic, and so we have shown the claim.

\end{proof}

\begin{theorem} 
\label{DensityLogSubanalytic}

Let $\pi:X\to Y$ be a morphism between affine algebraic varieties over $\RR$.
Let $U\subset X$ be a dense open set such that $\pi|_U$ is smooth. 
Suppose that $\w\in \G(U,\W_{U/Y})$ and that $f:X(\RR)\to \RR$ is a nonnegative continuous compactly-supported semialgebraic function.
Then on $Y$, the density function $y\mapsto \int_{\pi^{-1}(y)\cap U} f|\w|$ log-subanalytic away from the locus where it is infinite, which is a semialgebraic subset of $Y$.

\end{theorem}

\begin{proof}

We may assume that $Y\subset \AA^M$ and $X\subset Y\times \AA^N$ and that the morphism is the projection.

Up to decomposition and dilation, we may assume that $f$ is supported inside $(-1,1)^{M+N}$.
After further decomposition, we may assume there is a subset $I\subset\{1,\cdots,N\}$ such that for any $y\in Y$, the projection $\pi_I:\pi^{-1}(y)\cap U\cap (-1,1)^{M+N}\to (-1,1)^I$ is etale. 
Let $\nu$ be the standard volume form on $B(-1,1)^I$. 
For any $y\in Y(\RR)$, we have the following:

$$\int_{(\pi^{-1}(y)\cap U\cap (-1,1)^{M+N})(\RR)} |\w| = \int_{z\in (-1,1)^I} \left( \sum_{x\in (\pi_I^{-1}(z)\cap\pi^{-1}(y)\cap U \cap (-1,1)^{M+N})} \left| \frac{\w}{\pi_I^*\nu}(x) \right| \right) |\nu| $$

The integrand on the right is semialgebraic, and via the embedding of $Y\subset \AA^M$, we are in a position to apply Theorem \ref{CluckersMillerTheorem} of Cluckers and Miller after checking that our parameterized integrand is measurable and integrable at each value of $y\in Y$. 
By Lemma \ref{InfiniteLocusExcision}, our integrand is measurable and the locus where it is not integrable is semialgebraic and may be removed from the domain.
Thus the density function $y\mapsto \int_{\pi^{-1}(y)\cap U} f|\w|$ from $Y(\RR) \setminus \infty(Y(\RR),f|\w|) \to \RR$ is log-subanalytic.

\end{proof}

\section{Reduction to a Curve}

In this section, we show that in order to prove our main theorem, it is enough to prove the following theorem:

\begin{theorem}[Modified Main Theorem]
\label{ModifiedMainTheorem}

Let $\cX$ be a Stein complex-analytic variety with a complex conjugation $\sigma$.
Let $\psi:\cX\to \CC$ be a (FRS) map of complex-analytic varieties which intertwines $\sigma$ with the usual complex conjugation on $\CC$. 
Assume that $\cX$ is Gorenstein.
Let $\w_\cX$ be a regular nowhere-vanishing $\sigma$-invariant holomorphic top differential form on the smooth locus $X^{sm}$ of $X$ and let $f:\cX^\sigma\to \RR$ be a nonnegative continuous compactly-supported function.
Denote the smooth locus of $\psi$ by $\cX^S$. 
Then the measure $\psi_*(f|\w_\cX|)$ on $\cX^\sigma$ has continuous density with respect to the measure $|dz|=dx$ on $\CC^\sigma=\RR$, which is given by 

$$ (\psi_*f)(y) := \int_{(\psi^{-1}(y)\cap \cX^S)^\sigma} f\cdot\left| \frac{\w_\cX}{\psi^*dz}|_{\psi^{-1}(y)\cap \cX^S} \right|.$$


\end{theorem}

See Appendix A for definitions and notions from complex-analytic varieties.

Our strategy will be as follows: in subsection one, we will reduce the proof of the main theorem to the case of $X,Y$ affine and (CSCD) measures on $X$ of the form $f|\w|$ for $f$ semialgebraic.
In subsection two, we will show that if our main theorem fails for affine $X,Y$ and semialgebraic $f$, then we can find a complex-analytic variety $\cX$ and map $\psi:\cX\to\CC$ satisfying the conditions of Theorem \ref{ModifiedMainTheorem} where $\psi_*(f|\w_\cX|)$ does not have continuous density.

\subsection{Reduction to \texorpdfstring{$X,Y$}{X,Y} Affine and \texorpdfstring{$f$}{f} Semialgebraic}

\begin{theorem}
\label{AffineReduction}

In proving Theorem \ref{MainTheorem}, it suffices to consider the scenario when $X$ and $Y$ are both affine.

\end{theorem}

\begin{proof}

Suppose $\varphi:X\to Y$ is an (FRS) morphism of $\RR$-varieties with $Y$ smooth and $X$ Gorenstein.
Let $m_X$ be a (CSCD) measure on $X$.

First we prove that it is enough to consider $Y$ affine.
Let $V_i$ be a finite affine open covering of $Y$, and let $X_i=X\times_Y V_i$. 
The $X_i$ form a finite open cover of $X$, and the $X_i(\RR)$ form a finite open cover of $X(\RR)$.
Construct on $X(\RR)$ a locally finite partition of unity $\Psi_j$ such that each $\Psi_j$ is compactly supported and the support of each $\Psi_j$ is contained in some $X_i(\RR)$.
Then $$\varphi_*m_X=\varphi_*\left(\sum_j \Psi_jm_X\right)=\sum_i(\varphi|_{X_i})_*\left(\sum_j\Psi_jm_X|_{X_i}\right)$$ and because only finitely many $\Psi_j$ are nonzero on the support of $m_X$, we may reduce from checking if $\varphi_*m_X$ is (CSCD) to checking that each $(\varphi|_{U_i})_*(\Psi_jm_X|_{U_i})$ is (CSCD).

Next we prove that it is enough to consider $X$ affine.
For any $X$, we may take a finite covering of $X$ by affine opens $U_i$.
The $U_i(\RR)$ form an open cover of $X(\RR)$.
Construct on $X(\RR)$ a locally finite partition of unity $\rho_j$ such that each $\rho_j$ is compactly supported and the support of each $\rho_j$ is contained in some $U_i(\RR)$.
We may now write $$\varphi_*m_X=\varphi_*\left(\sum_j \rho_jm_X\right)=\sum_i(\varphi|_{U_i})_*\left(\sum_j\rho_jm_X|_{U_i}\right)$$ and because only finitely many $\rho_j$ are nonzero on the support of $m_X$, we may reduce from checking if $\varphi_*m_X$ is (CSCD) to checking that each $(\varphi|_{U_i})_*(\rho_jm_X|_{U_i})$ is (CSCD).

Since each of these reduction steps maintains the setting of a Gorestein source, smooth target, and (FRS) map, we have reduced our main theorem to the case of a map between affine $\RR$-varieties.

\end{proof}

Suppose $\varphi:X\to Y$ is (FRS) with $X$ Gorenstein and $Y$ smooth.
Now that we have reduced to the case of affine varieties, we may assume that $m_X$ is given by $f|\w_X|$ for $f:X(\RR)\to\RR$ a compactly-supported continuous function and $\w_X\in\W_X(X)$.
In order to prove that $\varphi_*m_X$ is (CSCD), we will show that the density of $f|\w_X|$ with respect to $|\w_Y|$ is continuous, where $\w_Y\in\W_Y(Y)$ is a nonvanishing top differential form on $Y$.

Let $X^S$ be the smooth locus of $\varphi$.
Since all fibers of $\varphi$ are reduced and $\varphi$ is flat, the smooth locus of $\varphi^{-1}(y)$ is equal to $X^S\cap \varphi^{-1}(y)$. 
Since the restriction of $f$ to $\varphi^{-1}(y)(\RR)$ is nonnegative, continuous, and compactly supported, Lemma \ref{RationalSingFinIntegral} implies that each integral $\int_{X^S\cap\varphi^{-1}(y)(\RR)}f\left|\frac{\w_X}{\varphi^*\w_Y}\right|$ is convergent. 
We denote the function $y\mapsto \int_{X^S\cap\varphi^{-1}(y)(\RR)}f\left|\frac{\w_X}{\varphi^*\w_Y}\right|$ by $\varphi_*f$ (this depends on $\w_X$ and $\w_Y$, despite their failure to appear in the notation). 
By Lemma \ref{AbsolutelyContinuousUpgraded}, $\varphi_*f$ is a function representing the density of $\varphi_*(f|\w_X|)$ with respect to $|\w_Y|$.
In order to show that $\varphi_*m_X$ is (CSCD), it suffices to show that $\varphi_*f$ is continuous.

The next theorem establishes that it is enough to consider the case of $f$ semialgebraic and nonnegative.

\begin{theorem}
\label{ReductionToSemialgebraic}

Suppose $X\subset\AA^n_\RR,Y\subset\AA^m_\RR$ are affine and $\varphi:X\to Y$ is (FRS).
Fix top differential forms $\w_X\in\W_X(X)$ and $\w_Y\in\W_Y(Y)$ and suppose $\w_Y$ is nowhere-vanishing. 
If the density of $\varphi_*(g|\w_X|)$ with respect to $|\w_Y|$ is continuous for every nonnegative, continuous, compactly supported, semialgebraic function $g:X(\RR)\to\RR$, then the density of $\varphi_*(f|\w_X|)$ is continuous for all continuous compactly supported functions $f:X(\RR)\to\RR$.

\end{theorem}

\begin{proof}

First we note that it is enough to check the continuity of $\varphi_*f$ when $f$ is nonnegative, as any continuous compactly supported function may be written as the difference of two continuous compactly supported nonnegative functions.
If $f$ is nonnegative, we may pick a sequence of nonnegative continuous compactly supported semialgebraic functions $f_i:X(\RR)\to\RR$ converging uniformly to $f$ with the additional property that $f(x)\geq f_i(x)$ for all $x\in X(\RR)$ and all $i$ by Corollary \ref{SWPunchline}.
We may also fix $b:X(\RR)\to\RR$ a continuous, compactly supported, semialgebraic bump function taking the value $1$ on $\operatorname{Supp} f$ by Lemma \ref{SemialgebraicCompactTrick}.
It is enough to prove that the density of $\varphi_*f_i$ converges uniformly to the density of $\varphi_*f$.

In order to show that $\varphi_*f_n$ converges uniformly to $\varphi_*f$ on $Y(\RR)$, note that $$\left|\int_{(\varphi^{-1}(y)\cap X^{s})(\RR)} f_n\left|\frac{\w_X}{\varphi^*\w_Y}\right| - \int_{(\varphi^{-1}(y)\cap X^{s})(\RR)} f\left|\frac{\w_X}{\varphi^*\w_Y}\right|\right| \leq \int_{(\varphi^{-1}(y)\cap X^{s})(\RR)} |f-f_n| \left|\frac{\w_X}{\varphi^*\w_Y}\right|. $$
Combining the above observation with the statement that for all $\e>0$ we may pick a $N$ such that for all $n\geq N$ and all $x\in X(\RR)$, we have $|f(x)-f_n(x)|\leq \e$, we see that the right hand side of the above equation is in fact bounded above by $\e \int_{(\varphi^{-1}(y)\cap X^{s})(\RR)} b \left|\frac{\w_X}{\varphi^*\w_Y}\right|$.
But this is exactly the density of $\e(\varphi_*b)$, which is continuous by our assumption.
Since $\varphi_*b$ is continuous and compactly supported, it has a finite maximum, and therefore the density of $\e\varphi_*b$ uniformly converges to zero as $\e$ does.
This means that the density of $\varphi_*f_n$ converges uniformly to density of $\varphi_*f$, and therefore the latter is continuous.

\end{proof}

\subsection{Construction of the Scenario of Theorem \ref{ModifiedMainTheorem}}

By results of the previous section, it is enough to show that for $\varphi:X\to Y$ an (FRS) map of affine $\RR$-varieties with $X\subset\AA^N$ Gorenstein and $Y\subset\AA^M$ smooth, top differential forms $\w_X\in\W_X(X)$ and $\w_Y\in\W_Y(Y)$, and nonnegative continuous compactly supported semialgebraic $f:X(\RR)\to\RR$ that $\varphi_*f:y\mapsto \int_{X^S\cap\varphi^{-1}(y)(\RR)}f\left|\frac{\w_X}{\varphi^*\w_Y}\right|$ is continuous.
In this section, we explain how to produce the scenario of Theorem \ref{ModifiedMainTheorem} given the above data.
For the rest of the section, we will assume that $X,Y,\varphi,\w_X,\w_Y,$ and $f$ are all as above.


Assume that $\varphi_*f$ is not continuous at a point $y\in Y(\RR)$.
Since $\varphi$ is (FRS), $\infty(X,f|\w|)$ is empty by an application of Lemma \ref{RationalSingFinIntegral}.
As $\infty(X,f|\w|)$ is empty, we may apply Theorem \ref{DensityLogSubanalytic} to obtain that  $\varphi_*f$ is log-subanalytic on $Y(\RR)\subset \AA^N_\RR(\RR)$. 

By Theorem \ref{LogSubanalyticDiscontinuity}, we may find an analytic map $g:\RR\to Y(\RR)$ such that $g(\RR)$ contains $y$ and $\varphi_*f|_{g(\RR)}$ is discontinuous at $y$.
Since $g: \RR \to Y(\RR)$ is a real-analytic map of smooth real-analytic manifolds, we may naturally upgrade it to a complex-analytic map $g^\CC$ from an open complex neighborhood of $\RR$ in $\CC$, which we denote $\Op(\RR)$, to $Y(\CC)$.
Without loss of generality, we may assume that $\Op(\RR)$ has a complex conjugation fixing $\RR$.
Recall that $\varphi^\RR:X(\RR)\to Y(\RR)$ is induced by $\varphi:X\to Y$, which also induces a map $\varphi^\CC:X(\CC)\to Y(\CC)$. 

Let $\cX=X(\CC)\times_{Y(\CC)} \Op(\RR)$, where we form the fiber product $X(\CC)\times_{Y(\CC)} \Op(\RR)$ using the maps $g^\CC: \Op(\RR) \to Y(\CC)$ and $\varphi^\CC:X(\CC)\to Y(\CC)$.
$X(\CC)\times_{Y(\CC)} \Op(\RR)$ is Stein by Corollary \ref{SteinConsequences}.
Composing the projection $\pi: X(\CC)\times_{Y(\CC)} \Op(\RR)\to \Op(\RR)$ with the inclusion $\Op(\RR)\hookrightarrow\CC$, we obtain $\psi:X(\CC)\times_{Y(\CC)} \Op(\RR)\to \CC$.
Since each of $X(\CC),Y(\CC),\Op(\RR)$ have complex conjugations and the maps $\varphi^\CC,g^\CC$ are intertwiners for these conjugations, the fiber product $\cX=X(\CC)\times_{Y(\CC)} \Op(\RR)$ also has a complex conjugation.
We also note that $\psi:\cX\to\CC$ is an intertwiner of the complex conjugations on $\cX$ and $\CC$.
Finally, the complex conjugation on $X(\CC)\times_{Y(\CC)} \Op(\RR)$ fixes exactly the real-analytic fiber product $X(\RR)\times_{Y(\RR)} \RR\subset X(\CC)\times_{Y(\CC)} \Op(\RR)$.

To produce a compactly supported, continuous, nonnegative $f:\cX^\sigma=X(\RR)\times_{Y(\RR)} \RR\subset X(\CC)\times_{Y(\CC)} \Op(\RR)\to \RR$, we may restrict $f:X(\RR)\to\RR$ to $X(\RR)\times_{Y(\RR)} \RR\subset X(\RR)$.
In order to specify $\w_\cX$, recall that $\frac{\w_X}{\varphi^*\w_Y}$ is a section of $\W_{X/Y}$ over the smooth locus of $\varphi$.
This naturally produces a section of the bundle of relative differential forms for the map $\psi:\cX\to\CC$, and we also call this section $\frac{\w_X}{\varphi^*\w_Y}$.
We may then specify $\w_\cX=(\frac{\w_X}{\varphi^*\w_Y})\wedge \psi^*dz$, which is a regular holomorphic conjugation-invariant top differential form on $\cX^{S}$.
This choice of $\w_\cX$ implies that the density function on $X(\RR)\times_{Y(\RR)} \RR$ is the same as the density function along $g(\RR)\subset Y(\RR)$. 
In particular, this density function $\psi_*f$ is discontinuous if and only if the original density function $\varphi_*f$ was.

To complete the reduction it remains to show that $\pi:X(\CC)\times_{Y(\CC)} \Op(\RR)\to \Op(\RR)$ is (FRS) and $X(\CC)\times_{Y(\CC)} \Op(\RR)$ is Gorenstein.

\begin{lemma}
\label{FiberProductFRS}

$\pi:X(\CC)\times_{Y(\CC)} \Op(\RR)\to \Op(\RR)$ is (FRS).

\end{lemma}

\begin{proof}

Consider the following diagram of completed local rings:

\begin{center}
\begin{tikzpicture}[node distance=3cm,auto]
  \node (XL0) {$\hat{\cO}_{X(\CC)\times_{Y(\CC)} \Op(\RR), p}^{an}$};
  \node (X) [right of=XL0]{$\hat{\cO}_{X(\CC),x}^{an}$};
	\node (L0) [below of=XL0] {$\hat{\cO}_{\Op(\RR),r}^{an}$};
  \node (Y) [right of=L0]{$\hat{\cO}_{Y(\CC),y}^{an}$};
	\draw[<-] (XL0) to node {}(X);
	\draw[<-] (XL0) to node {$g$}(L0);
	\draw[<-] (X) to node {$\hat{f}$}(Y);
	\draw[<-] (L0) to node {} (Y);
\end{tikzpicture}
\end{center}

Note that the diagram is not cartesian, but the top left entry $\hat{\cO}_{X(\CC)\times_{Y(\CC)} \Op(\RR), r}^{an}$ is a completion and localization of $\hat{\cO}_{X(\CC),x}^{an}\otimes_{\hat{\cO}_{Y(\CC),y}^{an}} \hat{\cO}_{\Op(\RR),r}^{an}$, so it is flat over $\hat{\cO}_{\Op(\RR),r}^{an}$ as $\hat{\cO}_{X(\CC),x}^{an}$ is flat over $\hat{\cO}_{Y(\CC),y}^{an}$.

To show that the fibers are reduced and have rational singularities, note the fibers of $\pi$ as complex-analytic varieties are the analytifications of the complex points of the fibers of $f$ as algebraic varieties, and therefore each fiber is reduced and has rational singularities.

\end{proof}


\begin{lemma}
\label{ModifiedGorenstein}

$X(\CC)\times_{Y(\CC)} \Op(\RR)$ is Gorenstein.

\end{lemma}

\begin{proof}

It is enough to show that for each point $p\in X(\CC)\times_{Y(\CC)} \Op(\RR)$ the completion of the local ring $\hat{\cO}_{X(\CC)\times_{Y(\CC)} \Op(\RR), p}^{an}$ is Gorenstein. 
We recall the following classical result on Gorenstein local rings \cite{GorensteinResult}: 

Let $(A,\mathfrak{m})$ and $(B,\mathfrak{n})$ be local rings with $\psi:A\to B$ a local homomorphism making $B$ flat over $A$. 
Then $B$ Gorenstein is equivalent to $A$ and $B/\mathfrak{m}B$ Gorenstein. 

We will apply this lemma to the map $g: \hat{\cO}^{an}_{\RR,r} \to \hat{\cO}^{an}_{X(\CC)\times_{Y(\CC)} \Op(\RR),(x,r)}$ from the left column of the diagram in Lemma \ref{FiberProductFRS}.
In our case, flatness has already been verified in Lemma \ref{FiberProductFRS}, and $\hat{\cO}^{an}_{\RR,r}$ is Gorenstein as $\RR$ is smooth. 
It remains to show that $\hat{\cO}^{an}_{X(\CC)\times_{Y(\CC)} \Op(\RR),(x,r)}/\mathfrak{m}_{\hat{\cO}^{an}_{\Op(\RR),r}}\hat{\cO}^{an}_{X(\CC)\times_{Y(\CC)} \Op(\RR),(x,r)}$ is Gorenstein. 
But this ring is isomorphic to the local ring of $x$ in the fiber of $X(\CC)\times_{Y(\CC)} \Op(\RR)$ over $l$, i.e. $\hat{\cO}^{an}_{(X(\CC)\times_{Y(\CC)} \Op(\RR))_r,x}\cong \hat{\cO}^{an}_{X_{f(x)},x}\cong \hat{\cO}_{X_{f(x)},x}$, which is just the completion of the local ring of the point in the algebraic fiber over $x$.
Applying the statement above to $\cO_{Y,y}\stackrel{f^\#}{\to} \cO_{X,x}$, we see that this local ring is Gorenstein.
Since a local ring is Gorenstein if and only if its completion is, we are done.

\end{proof}

From the combination of Lemmas \ref{FiberProductFRS} and \ref{ModifiedGorenstein}, we may conclude that we are in the situation of Theorem \ref{ModifiedMainTheorem} and we have completed our reduction.

\section{Reduction to a Local Model}

In this section, we reduce from the case of Theorem \ref{ModifiedMainTheorem} to computing the pushforward of monomial measures under monomial maps $\CC^n\to\CC$.
In subsection one, we will prove that if $\psi:\cX\to\cY$ is an (FRS) map of complex-analytic varieties with $\cY$ smooth, then $\cX$ has rational singularities.
In subsection two, we will prove a regularity theorem for meromorphic differential forms on resolutions of singularities of complex-analytic varieties with rational singularities.
In subsection three, we combine these two results to complete the reduction.

\subsection{\texorpdfstring{$\psi:\cX\to\cY$}{Psi:X to Y} (FRS) and \texorpdfstring{$\cY$}{Y} smooth implies \texorpdfstring{$\cX$}{X} has rational singularities}

In the algebraic setting, Elkik's Theorem shows that if $f:X\to Y$ is flat with $Y$ smooth, then any rational singularity $x$ of $X_y$ is also a rational singularity of $X$ \cite{Elkik}.
This theorem enables us to conclude that if $f:X\to Y$ is an algebraic morphism which is (FRS) with $Y$ smooth, then $X$ has rational singularities.
We are interested in an analytic analogue of this theorem in order to perform later computations.
Our goal is to show that if $f:\cX\to \cY$ is a map of complex-analytic varieties with $\cY$ smooth and $f$ (FRS), then $\cX$ has rational singularities.
We will apply this to conclude that $\cX$ in Theorem \ref{ModifiedMainTheorem} has rational singularities.
The proof will follow the one in Elkik's paper \cite{Elkik} very closely, with some mild departures due to the analytic setting.

\begin{lemma}
\label{FlatNCM}

Let $f:\cX\to \cS$ be a flat holomorphic map of complex-analytic spaces with $\cS$ smooth.
Let $s\in \cS$ and $x\in \cX_s$.
If $\cX_s$ is normal and Cohen-Macaulay at $x$, then $\cX$ is normal and Cohen-Macaulay at $x$.

\end{lemma}

\begin{proof}

We apply a characterization for flatness from Hironaka \cite{HironakaFlat}.
Let $C(\cX,F)$ denote the normal cone to $F$ in $\cX$, where $F$ is a closed analytic subvariety of $\cX$.
If $f:\cX\to \cS$ is a holomorphic map between two complex-analytic spaces with $f(x)=s$, then $f$ is flat at $x$ if and only if in a sufficiently small neighborhood of $x$, we have $C(\cX,\cX_s)\cong \cX_s\times C(\cS,s)$.

Recall that a local ring is normal or Cohen-Macaulay if and only if its completion is.
To check if $\cX$ is normal and Cohen-Macaulay at $x$, it is enough to check at the completed local ring of $x$.
From the construction of the normal cone, we recover that $\widehat{\cO}^{an}_{\cX,x}\cong \widehat{\cO}^{an}_{C(\cX,X_s),x}$, where the completion is taken with respect to the maximal ideal.
Applying Hironaka's characterization to the map $f:\cX\to \cS$, we see that $C(\cX,X_s)\cong \cX_s\times C(\cS,s)$ in some small-enough open neighborhood of $x\subset \cX_s$, which implies that $\cO^{an}_{C(\cX,\cX_s),x}\cong \cO^{an}_{\cX_s,x}\otimes\cO^{an}_{C(\cS,s),s}$. 
Each of these local rings are normal and Cohen-Macaulay, so their tensor product is as well.
Thus $\widehat{\cO}^{an}_{C(\cX,X_s),x}$ and therefore $\widehat{\cO}^{an}_{\cX,x}$ are normal and Cohen-Macaulay.
This proves that that $\cX$ is normal and Cohen-Macaulay at $x$

\end{proof}

The following analytic analogue of the Grauert-Riemanschneider Vanishing Theorem will be useful:

\begin{lemma}
\label{AnalyticGR}

Suppose $\pi:\widetilde{\cX}\to \cX$ is a resoultion of singularities of $\cX$ a complex-analytic space.
Then $R^i\pi_*\W_{\widetilde{\cX}}=0$ for $i>0$, and $\pi_*\W_{\widetilde{\cX}}$ is independent of the resolution of singularities.

\end{lemma}
	
\begin{proof}

See Lemma A.2 in \cite{SilvaVanishing}, or Corollary 1 in \cite{Takegoshi} for a slightly more general version.

\end{proof}

\begin{theorem}
\label{AnalyticElkik}

Let $f:\cX\to \cS$ be a flat holomorphic map of complex-analytic spaces with $\cS$ smooth and $\cX$ reduced.
Let $s\in S$.
If $x\in \cX$ is a rational singularity of $\cX_s$, then $x$ is a rational singularity of $\cX$.

\end{theorem}

\begin{proof}

As each $\cX_s$ has rational singularities, Theorem \ref{MoreRationalSingularitiesAnalytic} gives that $\cX_s$ is normal and Cohen-Macaulay at $x\in X_s$.
By Lemma \ref{FlatNCM}, $\cX$ is normal and Cohen-Macaulay at $x$.
After restricting to a neighborhood of $x\in \cX$, we may assume that $\cX$ is normal and Cohen-Macaulay.

Since $\cX$ is Cohen-Macaulay, the dualizing complex $\W_\cX^\bullet$ is isomorphic to a complex concentrated in a single degree.
By abuse of notation, let $\W_\cX$ be a sheaf such that $\W_\cX^\bullet$ is isomorphic to the complex consisting of $\W_\cX$ in the correct degree and zero elsewhere.
Let $\pi:\widetilde{\cX}\to \cX$ be a resolution of singularities.
By Lemma \ref{AnalyticGR}, an analytic analogue of the Grauert-Riemanschneider Vanishing Theorem, we have that $R^i\pi_*\W_{\widetilde{\cX}}=0$ for $i>0$, and $\pi_*\W_{\widetilde{\cX}}$ is independent of the choice of resolution.

Denote $i: \cO_\cX\to R^\bullet\pi_*\cO_{\widetilde{\cX}}$ the natural map.
Dualizing, we obtain $j: R^\bullet\pi_*\W_{\widetilde{\cX}}[\dim \widetilde{\cX}]\to \W_\cX^\bullet$, which by the discussion in the previous paragraph, is a map between two complexes concentrated in degree $\dim \widetilde{\cX}$. 
We will refer to this map as an honest map of sheaves $j:\pi_*\W_{\widetilde{\cX}}\to \W_\cX$.
By Theorem \ref{MoreRationalSingularitiesAnalytic}, in order to prove that $x$ is a rational singularity of $\cX$, it suffices to prove that $j$ is surjective in a neighborhood of $x$.

We will proceed by induction on $\dim \cS$.
Suppose $\dim \cS=1$.
Let $s\in \cS$.
Suppose $V\subset \cS$ is an open neighborhood of $s$ with a choice of $t\in \G(V,\cO_\cS)$ such that $t$ generates the maximal ideal in the stalk of $\cO_\cS$ at $s$.
We obtain the following exact sequences on $\widetilde{\cX}$ and $\cX$, respectively:

$$0\to \W_{\widetilde{\cX}} \stackrel{(f\circ\pi)^*t}{\to} \W_{\widetilde{\cX}} \to \W_{\widetilde{\cX}_s} \to 0$$
$$0 \to \W_\cX \stackrel{f^*t}{\to} \W_\cX \to \W_{\cX_s} \to 0.$$

Taking direct images of the first exact sequence and applying Lemma \ref{AnalyticGR}, we see that it is still exact. 
Using the morphism $j: \pi_*\W_{\widetilde{\cX}}\to \W_\cX$, we obtain the following diagram:

\begin{center}
\begin{tikzpicture}[node distance=2cm,auto]
  \node (0tl) {$0$};
  \node (t1) [right of=0tl]{$\pi_*\W_{\widetilde{\cX}}$};
  \node (t2) [right of=t1]{$\pi_*\W_{\widetilde{\cX}}$};
	\node (t3) [right of=t2]{$\pi_*\W_{\widetilde{\cX}_s}$};
	\node (0tr) [right of=t3]{$0$};
	\node (0bl) [below of=0tl] {$0$};
  \node (b1) [right of=0bl]{$\W_{\cX}$};
  \node (b2) [right of=b1]{$\W_{\cX}$};
	\node (b3) [right of=b2]{$\W_{\cX_s}$};
	\node (0br) [right of=b3]{$0$};
  \draw[->] (0tl) to node {}(t1);
	\draw[->] (t1) to node {$f^*t$} (t2);
	\draw[->] (t2) to node {}(t3);
	\draw[->] (t3) to node {}(0tr);
	\draw[->] (0bl) to node {}(b1);
	\draw[->] (b1) to node {$f^*t$} (b2);
	\draw[->] (b2) to node {}(b3);
	\draw[->] (b3) to node {}(0br);
  \draw[->] (t1) to node {$j$}(b1);
	\draw[->] (t2) to node {$j$}(b2);
	\draw[->] (t3) to node {$k$}(b3);
\end{tikzpicture}
\end{center}

where $k$ is induced by $j$.

By taking a resolution of the pair $(\cX,\cX_s)$ we may assume that $\cX_s'$, the strict transform of $\cX_s$ in $\widetilde{\cX}$, is a resolution of singularities of $\cX_s$.
As $\cX_s$ was assumed to have a rational singularity at $x$, we have that $\pi_*\W_{\cX_s'}\to \W_{\cX_s}$ is surjective in a neighborhood of $x$.
Noting that $\pi_*\W_{\cX_s'}$ is naturally a subsheaf of $\pi_*\W_{\widetilde{\cX}_s}$, we see that $k$ is surjective in a neighborhood of $x$.
Applying the Snake Lemma to the above diagram of short exact sequences, we see that $\coker j\stackrel{f^*t}{\to}\coker j \to \coker k$ is an exact sequence.
Restricting to a neighborhood of $x$ where $k$ is surjective, we have that $\coker j\stackrel{f^*t}{\to}\coker j \to 0$ is exact.
By Nakayama's Lemma, $\coker j=0$ in this neighborhood, and therefore $j$ is surjective in our chosen neighborhood of $x$.
Therefore we have shown that $x$ is a rational singularity of $\cX$.

To complete the induction, suppose that the statement is true for $\dim \cS \leq n$. 
Let $\cS$ be smooth of dimension $n+1$.
Pick $V\subset \cS$ an open coordinate neighborhood of $s$ with coordinates $t_1,\cdots,t_{n+1}$ generating the maximal ideal of $\cO_\cS$ at $s$.
Let $V_{n+1}=\{t_1=\cdots=t_n=0\}\subset V$.
Then $\cX\times_s V_{n+1}\to V_{n+1}$ is a flat holomorphic map with $V_{n+1}$ smooth and of dimension one with the same fibers as $f:\cX\to \cS$.
Thus each point $x\in \cX_s$ which is a rational singularity is also a rational singularity of $\cX\times_s V_{n+1}$.

On the other hand, let $V'\subset\CC^n$ be obtained from $V$ by projecting away the final coordinate, i.e. $V'$ is the image of $V$ under the map $$V\stackrel{(t_1,\cdots,t_{n+1})}{\hookrightarrow} \CC^{n+1}\stackrel{proj_{t_{n+1}}}{\to} \CC^n.$$
Then the composite $f^{-1}(V)\to V'$ is flat with fiber over zero isomorphic to $\cX\times_s V_{n+1}$, and so by the induction hypothesis, $x$ is a rational singularity of $\cX$.

\end{proof}

\subsection{Regularity on the Exceptional Divisor}

The next theorem will be used to establish key bounds in the following section.

\begin{theorem}

\label{RSRegularEDivisor}

Let $\cY$ be a complex-analytic variety which is Stein, Gorenstein, and has rational singularities.
Let $\varphi:\cY\to \CC$ be a holomorphic function, and let $\cY_0=\varphi^{-1}(0)$.
Let $\pi:\widetilde{\cY}\to \cY$ be a resolution of singularities, let $\widetilde{\varphi}=\varphi\circ\pi$, and $\widetilde{\cY}_0=\widetilde{\varphi}^{-1}(0)$.
If $\widetilde{\eta}$ is a meromorphic top form defined on $\widetilde{\cY}$ so that $\widetilde{\varphi}\widetilde{\eta}$ is holomorphic, then $\widetilde{\eta}$ is holomorphic on $\widetilde{\cY}_0\setminus \cY_0'$.

\end{theorem}

\begin{proof}

Since $\widetilde{\varphi}\widetilde{\eta}$ is regular and $\cY$ has rational singularities, there exists a top differential form $\z$ on $\cY$ such that $\z$ agrees with $\widetilde{\varphi}\widetilde{\eta}$ on an open dense set. 
Let $\eta=\frac{\z}{\varphi}$. 
Let $\w_0$ be the residue of $\eta$ along $\cY_0$. 
Since $\cY_0$ has rational singularities and $\cY_0'$ is a resolution of singularities of $\cY_0$, we may find a top form $\widetilde{\w}_0$ on $\cY_0'$ which agrees with $\w_0$ along an open dense set.

In order to find a meromorphic top form $\widetilde{\w}$ on $\widetilde{\cY}$ so that the residue of $\widetilde{\w}$ along $\cY_0'$ is $\widetilde{\w}_0$, we need to show that after taking global sections over $\widetilde{\cY}$, the residue exact sequence $$0\to \W_{\widetilde{\cY}}\to \W_{\widetilde{\cY}}([\cY_0']) \to i_*\W_{\cY_0'} \to 0$$ is still exact. 
Thus it suffices to show that $H^1(\widetilde{\cY},\W_{\widetilde{\cY}})=0$.

Since $\cY$ has rational singularities, we have that $\pi_*\W_{\widetilde{\cY}}\cong \W_\cY$ and $R^p\pi_*\W_{\widetilde{\cY}}=0$ for $p>0$.
This implies that $H^i(\widetilde{\cY},\W_{\widetilde{\cY}})\cong H^i(\cY,\pi_*\W_{\widetilde{\cY}})$ for all $i$.
Since $\cY$ is Stein, we have that $H^i(\cY,\pi_*\W_{\widetilde{\cY}})=0$ for all $i>0$ by an application of Cartan's Theorem B (see Theorem \ref{CartanB}).
Thus we may choose a meromorphic top form $\widetilde{\w}$ on $\widetilde{\cY}$ so that the residue of $\widetilde{\w}$ along $\cY_0'$ is $\widetilde{\w}_0$ and $\widetilde{\w}$ is regular away from $\cY_0'$.

Now consider $\widetilde{\d}=\widetilde{\varphi}(\widetilde{\eta}-\widetilde{\w})$.
By construction, $\widetilde{\d}$ is holomorphic on $\widetilde{\cY}$, so we may choose a holmorphic differential form $\d$ on $\cY$ which agrees with $\widetilde{\d}$ on an open dense set.
Computing residues along $\cY^{sm}\cap \cY_0$, we see that $\d$ vanishes on $\cY^{sm}\cap \cY_0$, so $\d$ vanishes along all of $\cY_0$ and is therefore divisible by $\varphi$.
This implies that $\widetilde{\d}$ is divisible by $\widetilde{\varphi}$, or that $\widetilde{\eta}-\widetilde{\w}$ is holomorphic.
This shows that $\widetilde{\eta}$ has no poles along $\widetilde{\cY}\setminus \widetilde{\cY}_0'$.

\end{proof}

\subsection{Reduction}

Suppose we are in the situation of Theorem \ref{ModifiedMainTheorem}, with a Stein, Gorenstein complex-analytic variety $\cX$ equipped with a complex conjugation $\sigma$, a $\sigma$-invariant holomorphic differential form $\w_\cX\in\W_\cX(\cX^{sm})$, a compactly-supported continuous function $f:\cX^\s\to\RR$ and an (FRS) map $\psi:\cX\to\CC$ which intertwines $\sigma$ with the usual complex conjugation on $\CC$.
Denote $\cX_0=\psi^{-1}(0)$, and let $\widetilde{\cX}\stackrel{\pi}{\to} \cX$ be a complex-conjugation equivariant strong resolution of singularities (see Theorem \ref{HironakaAlgebraicResolution}) of the pair $(\cX,\cX_0)$, so that $\widetilde{\cX}$ inherits a complex conjugation. 
The singular values of $\widetilde{\psi}=\psi\circ\pi$ form an analytic closed subset of $\CC$ with no accumulation points, therefore we may pick an open neighborhood of the origin $U\subset\CC$ such that $0$ is the only singular value inside $U$. 
In this neighborhood, for any $t\in U$, $\widetilde{\psi}^{-1}(t)$ is smooth and the map $\widetilde{\psi}^{-1}(t)\to\psi^{-1}(t)$ is a resolution of singularities. 
By a standard argument, the strict transform $\cX_0'$ of $\cX_0$ is a resolution of singularities of $\cX_0$.

Let $W\subset \cX$ be an open dense subset where $\pi$ is an isomorphism. 
Since $\cX$ has rational singularities by Theorem \ref{AnalyticElkik}, we may choose a top differential form $\w_{\widetilde{\cX}}\in \W_{\widetilde{\cX}}(\widetilde{\cX})$ agreeing with $\w_\cX$ on $W$, i.e. $\w_{\widetilde{\cX}}|_{\pi^{-1}W}=\pi^*(\w_{\cX}|_W)$. 

Since $\widetilde{\cX}_0$ is a divisor with strict normal crossings inside a smooth variety, locally at every point $z$ in $\widetilde{\cX}_0$ there exists a coordinate system $x_1,\cdots,x_n$ such that both $\widetilde{\psi}$ and $\w_{\widetilde{\cX}}$ are monomial, i.e.

$$\widetilde{\psi} = \a x_1^{a_1}\cdots x_n^{a_n}$$
$$\w_{\widetilde{\cX}} = \b x_1^{b_1}\cdots x_n^{b_n}dx_1\wedge\cdots\wedge dx_n$$

where $\a,\b$ are units. 
Since $\cX_0'\cap\pi^{-1}(W)$ is isomorphic to the reduced variety $W\cap \cX_0$, we get that $\cX_0'$ is reduced by Serre's criteria for reducedness and the fact that $\cX_0'$ is Cohen-Macaulay. 
For points $z\in \cX_0'$, we have that $\cX_0'$ is locally the zero locus of one $x_i$, which we may assume to be $x_1$, implying that $a_1=1$. 
By Theorem \ref{RSRegularEDivisor} applied to $\frac{\w_{\widetilde{\cX}}}{\widetilde{\psi}}$, if $z\in \cX_0'$, we have $a_i\leq b_i$ for $i\geq 2$ and if $z\notin \cX_0'$, then $a_i\leq b_i$ for all $i$.

Setting $\widetilde{f}=f\circ\pi$ and computing measures, we have that $\psi_*(f|\w_\cX|)=\widetilde{\psi}_*(\widetilde{f}|\w_{\widetilde{\cX}}|)$.
This means that in order to compute $\psi_*(f|\w_\cX|)$, it is enough to compute $\widetilde{\psi}_*(\widetilde{f}|\w_{\widetilde{\cX}}|)$.
The following partition of unity argument reduces the computation of $\widetilde{\psi}_*(\widetilde{f}|\w_{\widetilde{\cX}}|)$ to the case where $\widetilde{\cX}=\CC^n$ and both $\widetilde{\psi}$ and $\w_{\widetilde{\cX}}$ are monomial.

For each $x\in\widetilde{\cX}$, choose a relatively compact connected open neighborhood $U_x\subset\widetilde{\cX}$ in the standard topology such that both $\widetilde{\varphi}$ and $\w_{\widetilde{\cX}}$ are monomial in this neighborhood.
Up to possibly shrinking $U_x$, we may assume that each of these is an analytic coordinate patch where $\a$ and $\b$ extend to units on the closure of $U_x$.
Since $\widetilde{\cX}$ is a locally compact second-countable Hausdorff space, there exists a countable, locally finite refinement of the covering $\{U_x\}_{x\in\cX}$, which we denote $\{U_i\}_{i\in I}$.
Since the $U_i$ cover $\widetilde{\cX}$, the sets $U_i^\s$ cover $\widetilde{\cX}^\s$, and we may form a partition of unity $\Phi_i$ subordinate to the cover $\{U_i^\s\}_{i\in I}$.

Let $A:=\operatorname{Supp} \widetilde{f}\subset \widetilde{\cX}^\sigma$.
Since $\pi$ is proper and $f$ is compactly supported, $A$ is compact.
As $\{U_i^\s\}_{i\in I}$ forms an open covering of $\widetilde{\cX}^\s$, we may select a finite collection $\{U_i^\s\}_{i\in I'}$ which forms an open cover of $A$.
Let $f_i=\Phi_if$.
Then $f|\w_{\widetilde{\cX}}|=\sum_{i\in I'} f_i|\w_{\widetilde{\cX}}|$, so it is enough to demonstrate that $\widetilde{\psi}_*(f_i|\w_{\widetilde{\cX}}|)$ has continuous density.
By the choice of $U_i$, this is the same as showing the pushforward of $f_i|\beta x_1^{b_1}\cdots x_n^{b_n}dx_1\cdots dx_n|$ along the map $\a x_1^{a_1}\cdots x_n^{a_n}:\CC^n\to\CC$ has continuous density.



\section{Proof for a Local Model}

In this section, we prove the following theorem:

\begin{theorem}
\label{MainTheoremLocal}

Denote the usual measure on $\RR^n$ by $dx^{\boxtimes n}$ and the usual measure on $\RR$ by $dx$.
Assume that $f:\RR^n\to\RR$ is a continuous compactly supported nonnegative function, $\a,\b:\RR^n\to\RR$ are analytic functions which are units on the support of $f$, and $A=(a_1,\cdots,a_n)$ and $B=(b_1,\cdots,b_n)$ are sequences of nonnegative integers with not all $a_i$ equal to zero.
Suppose either of the following two cases hold:

\begin{enumerate}
\item $a_i\leq b_i$ for all $i$, or
\item $a_1=1$ and $a_i\leq b_i$ for $i\geq 2$.
\end{enumerate}

Then the pushforward of the measure $f|\beta x^B|dx^{\boxtimes n}=f|\beta x_1^{b_1}\cdots x_n^{b_n}dx_1\cdots dx_n|$ along the map $\a x^A = \a x_1^{a_1}\cdots x_n^{a_n}:\RR^n\to\RR$ has continuous density with respect to the measure $dx$.

\end{theorem}

The following lemma reduces our task to a computationally simpler one:

\begin{lemma}
\label{ReductionToStandard}

In proving Theorem \ref{MainTheoremLocal}, it suffices to treat the case where $\a=1$, $\b=1$, and $\operatorname{Supp} f\subset [-1,1]^n$.

\end{lemma}

\begin{proof}

To prove that it suffices to consider $\b=1$, note that $f|\b x^B|=(f|\b|)|x^B|$ and $f|\b|$ is a continuous compactly supported nonnegative function.
To prove that it suffices to consider $\a=1$, we note that we may make an analytic change of coordinates on $\RR^n$ which changes the map $\a x^A$ to $x^A$.
This sends the measure $f|x^B|dx^{\boxtimes n}$ to $g|x^B|dx^{\boxtimes n}$, where $g$ is also continuous, compactly supported, and nonnegative.
Since $\operatorname{Supp} f$ is compactly supported, it is bounded, and so up to a dilation action which multiplies the measure by a positive constant, we may assume that $\operatorname{Supp} f$ is contained in $[-1,1]^n$.

\end{proof}

Next, we show that for continuous nonnegative $f$ with $\operatorname{Supp} f\subset [-1,1]^n$, the measure $(x^A)_*(f|x^B|dx^{\boxtimes n})$ has continuous density on $\RR^n\setminus \{0\}$.

\begin{lemma}
\label{GenericContinuous}

Assume that $f:\RR^n\to\RR$ is a continuous nonnegative function with $\operatorname{Supp} f\subset [-1,1]^n$, and $A=(a_1,\cdots,a_n)$ and $B=(b_1,\cdots,b_n)$ are sequences of nonnegative integers with not all $a_i$ equal to zero.

Then the pushforward of the measure $f|\beta x^B|dx^{\boxtimes n}=f|\beta x_1^{b_1}\cdots x_n^{b_n}dx_1\cdots dx_n|$ along the map $\a x^A = \a x_1^{a_1}\cdots x_n^{a_n}:\RR^n\to\RR$ has continuous density with respect to the measure $dx$ on the set $\RR\setminus\{0\}$.

\end{lemma}

\begin{proof}

Let $g_m$ be a continuous function taking the value $1$ on $(|x^A|)^{-1}([\frac{1}{m+1},1])$, vanishing outside $(|\a x^A|)^{-1}([\frac{1}{m+2},2])$, and so that $0\leq g_m\leq 1$ for all $m$.
Then $(x^A)_*(fg_m|x^B|)=(x^A)_*(fg_{m+1}|x^B|)=(x^A)_*(f|x^B|)$ on $(-m,-\frac{1}{m})\cup(\frac{1}{m},m)$, and by Lemma \ref{WellBehavedSmoothPushforward}, $(x^A)_*(fg_m|x^B|)$ is continuous for each $m$.

\end{proof}

It remains to show that the density is continuous at zero.
In order to do this, we introduce an auxiliary definition and recall a basic fact from measure theory.

\begin{definition}
\label{KBox}

Suppose $k$ and $n$ are positive integers.
Suppose $r_1,\cdots,r_n$ are integers satisfying $0\leq r_i<k$.
Define $\square_{k,r_1,\cdots,r_n}=[-1+\frac{2r_1}{2k+1},-1+\frac{2(r_1+1)}{2k+1}]\times\cdots\times[-1+\frac{2r_n}{2k+1},-1+\frac{2(r_n+1)}{2k+1}]$.

\end{definition}

\begin{lemma}
\label{BoxApprox}

For any continuous function $f$ supported inside $[-1,1]^n$ and any $\e>0$, there exists an integer $k>0$ and a function $f^\square$ which is linear combination of characteristic functions of sets of the form $\square_{k,r_1,\cdots,r_n}$ so that $$|f-f^\square|< \e1_{[-1,1]^n}$$ almost everywhere.
Additionally, the locus where the approximation fails is a union of perpendicular translates of coordinate hyperplanes by odd multiplies of $\frac{1}{2k+1}$.

\end{lemma}

\begin{proof}

This is just the definition of integrability.

\end{proof}

\begin{lemma}
\label{ReductionToBox}

Suppose $A,B$ are as in the situation of Theorem \ref{MainTheoremLocal}.
Suppose $(x^A)_*(1_{\square_{k,r_1,\cdots,r_n}}|x^B|dx^{\boxtimes n})$ is continuous in a neighborhood of $0$.

Then $(x^A)_*(f|x^B|dx^{\boxtimes n})$ is continuous with respect to $dx$.

\end{lemma}

\begin{proof}

Recall that the function $(x^A)_*(f|x^B|)(y)$ given by $y\mapsto\int_{(x^A)^{-1}(y)^S} f|x^B|\frac{dx_1\cdots dx_n}{dx}$ represents the density of the pushforward $(x^A)_*(f|x^B|)dx^{\boxtimes n}$ with respect to $dx$.
By Lemma \ref{GenericContinuous}, it is enough to show that $(x^A)_*(f|x^B|)(y)$ is continuous at $0$. 

Let $\e>0$.
We will show directly that $$ |(x^A)_*(f|x^B|)(0)-(x^A)_*(f|x^B|)(q) | <\e $$ for all $q$ in some small neighborhood of $0$.

According to Lemma \ref{BoxApprox}, we may choose a function $f^\square$ which is a linear combination of characteristic functions of boxes so that $f^\square$ approximates $f$ uniformly within any $\e_1>0$ except on a union of translates of coordinate hyperplanes.
By restricting to small enough $q\in\RR$, we may assume that the intersection of the set of non-uniform approximation with $(x^A)^{-1}(q)^S$ is of measure zero inside $(x^A)^{-1}(q)^S$ for each $q$.
Therefore $(x^A)_*(|f-f^\square|\cdot|x^B|)\leq \e_1(x^A)_*(1_{[-1,1]^n}|x^B|)$ on a small enough open neighborhood of $0\in\RR$.

Writing $$ |(x^A)_*(f|x^B|)(0)-(x^A)_*(f|x^B|)(q) | \leq \left|(x^A)_*(f|x^B|)(0)-(x^A)_*(f^\square|x^B|)(0) \right| + $$ $$ + \left|(x^A)_*(f^\square|x^B|)(0) -(x^A)_*(f^\square|x^B|)(q) \right| + $$ $$ + \left| (x^A)_*(f|x^B|)(q)-(x^A)_*(f^\square|x^B|)(q)  \right|, $$ it is enough to show that each term on the right hand side is less than $\frac{\e}{3}$ in a small-enough open neighborhood of zero.

Rewriting the first term as $(x^A)_*(|f-f^\square|\cdot|x^B|)(0)$, we see that this is less than $\e_1(x^A)_*(1_{[-1,1]^n}|x^B|)(0)$.
Since $(x^A)_*(1_{[-1,1]^n}|x^B|)$ is continuous in a neighborhood of $0$, it is bounded near $0$, and by fixing $q$ in some small open neighborhood of $0$ and picking $\e_1$ small we may make the first term less than $\frac{\e}{3}$.
The same argument applies to the third term, and as the second term is a linear combination of continuous functions, we may also bound it by $\frac{\e}{3}$ on a sufficiently small neighborhood of zero.

\end{proof}

The following lemma enables us to reduce checking the continuity of $(x^A)_*(1_{\square_{k,x_1,\cdots,x_n}}|x^B|dx^{\boxtimes n})$ in a neighborhood of $0$ to checking the continuity of $(x^A)_*(1_{[-1,1]^n}|x^B|dx^{\boxtimes n})$ in a neighborhood of $0$.

\begin{lemma}
\label{ReductionToBigBox}

Suppose $A$ and $B$ are as in the situation of Theorem \ref{MainTheoremLocal}.
If $(x^A)_*(1_{[-1,1]^n}|x^B|dx^{\boxtimes n})$ is continuous in a neighborhood of $0$, then $(x^A)_*(1_{\square_{k,x_1,\cdots,x_n}}|x^B|dx^{\boxtimes n})$ is continuous in a neighborhood of $0$.

\end{lemma}

\begin{proof}

Let $\square$ be any box of the form $[-1+\frac{2r_1}{2k+1},-1+\frac{2(r_1+1)}{2k+1}]\times\cdots\times[-1+\frac{2r_n}{2k+1},-1+\frac{2(r_n+1)}{2k+1}]$ for each $0\leq r_i<2k+1$ an integer.
Let $\widehat{\square}$ be the union of all reflections of $\square$ through all coordinate subspaces of $\RR^n$.
We may write $1_{\widehat{\square}}$ as a linear combination of characteristic functions of $n$-dimensional boxes of the form $SB_{k,s_1,\cdots,s_n}=[-\frac{2s_1+1}{2k+1},\frac{2s+1}{2k+1}]\times\cdots\times[-\frac{2s_n+1}{2k+1},\frac{2s_n+1}{2k+1}]$.

First, there exists a positive integer $N$ so that $N(x^A)_*(1_\square|x^B|)=(x^A)_*(1_{\widehat{\square}}|x^B|)$.
Now, writing $1_{\widehat{\square}}$ as a linear combination described above, we see that $(x^A)_*(1_\square|x^B|)$ is a linear combination of functions of the form $(x^A)_*(1_{SB_{k,s_1,\cdots,s_n}}|x^B|)$.
But each of these functions is equal to a multiple of $(x^A)_*(1_{[-1,1]^n}|x^B|)$ composed with a scaling of the input by the change of variables formula.
So $(x^A)_*(1_\square|x^B|)$ is continuous in a neighborhood of $0$.

\end{proof}

As the $f^\square$ picked in Lemma \ref{ReductionToBox} was a linear combination of functions of the form $1_{\square_{k,x_1,\cdots,x_n}}$, it is enough to analyze $(x^A)_*(1_{[-1,1]^n}|x^B|dx^{\boxtimes n})$.
Finally, we reduce from the case of $[-1,1]^n$ to the case of $[0,1]^n$.

\begin{lemma}
\label{ReductionToSmallBox}

Assume $A$ and $B$ are as in Theorem \ref{MainTheoremLocal}.
In order to show that $(x^A)_*(1_{[-1,1]^n}|x^B|dx^{\boxtimes n})$ has continuous density with respect to $dx$ in a neighborhood of zero, it suffices to show that the function given by $\frac{-d}{dq} \int_{(x^A)^{-1}(q,1)} 1_{[0,1]^n}x^B$ for $q\in(0,1)$ has continuous extension to $q=0$, which is zero if all entries of $A$ are even.

\end{lemma}

\begin{proof}

If there is an odd entry in $A$, then $(x^A)_*(1_{[-1,1]^n}|x^B|)$ at a point $q$ is the same as the density at a point $-q$, and the density for $q\in (0,1)$ is given by $\lim_{\e\to0} \frac{2^{n-1}}{2\e}\int_{(x^A)^{-1}(q-\e,q+\e)} 1_{[0,1]^n}x^B$.
If $A$ has only even entries, then $(x^A)_*(1_{[-1,1]^n}|x^B|)$ is zero for $q<0$, and the density for $q\in(0,1)$ is given by $\lim_{\e\to0} \frac{2^n}{2\e}\int_{(x^A)^{-1}(q-\e,q+\e)} 1_{[0,1]^n}x^B$.
To show that $(x^A)_*(1_{[-1,1]^n}|x^B|)$ has continuous density in a neighborhood of $0$, we may show that the function given by $\lim_{\e\to0} \frac{1}{2\e}\int_{(x^A)^{-1}(q-\e,q+\e)} 1_{[0,1]^n}x^B$ for $q\in(0,1)$ has a continuous extension to $q=0$, and that this extension is $0$ when $A$ has only even entries.
As $\lim_{\e\to0} \frac{1}{2\e}\int_{(x^A)^{-1}(q-\e,q+\e)} 1_{[0,1]^n}x^B=-\frac{d}{dq} \int_{(x^A)^{-1}(q,1)} 1_{[0,1]^n}x^B$, we have proven the lemma.

\end{proof}

We start by computing $\int_{(x^A)^{-1}(q,1)} 1_{[0,1]^n}x^B$ for $q\in(0,1)$.
First, we note that we may generalize to the case where $A$ and $B$ are sequences of nonnegative real numbers with not all $a_i=0$.
We note that we may immediately integrate away all $a_i$ which are zero at the price of a constant, so we may assume in the following that all $a_i$ are positive real numbers.

In the following, if $S=(s_1,\cdots,s_k)$ is a sequence of real numbers, we say that $S>c$ or $S\geq c$ for $c\in\RR$ if $s_i>c$ or $s_i\geq c$ for all $i$.

\begin{lemma}
\label{VolumeComputation}

Suppose $A=(a_1,\cdots,a_n)$ is a sequence of positive real numbers and $B=(b_1,\cdots,b_n)$ is a sequence of nonnegative real numbers.
Let $C=(c_1,\cdots,c_n)$ where $c_i=\frac{b_i+1}{a_i}$.
Let $q\in (0,1)$.
Let $h_k(z_1,\cdots,z_l)$ be the complete symmetric homogeneous polynomial in the variables $z_1,\cdots,z_l$ with the convention that $h_k(z_1,\cdots,z_l)=0$ for $k<0$ and $h_0(z_1,\cdots,z_l)=1$.
Let $V(A,B,q)$ denote the volume of $(x^A)^{-1}(q,1)$ with respect to the measure $1_{[0,1]^n}x^Bdx^{\boxtimes n}$.
Then $V(A,B,q)$ is continuous in $A$, $B$, and $q$ and may be written in any of the following forms:

\begin{enumerate}
\item $V(A,B,q)=\frac{1}{\prod_{i=1}^n a_i}\displaystyle\sum_{i=0}^\infty \frac{h_{i-n}(c_1,\cdots,c_n)(\log q)^i}{i!}$,
\item $V(A,B,q)=\frac{1}{\prod_{i=1}^n a_i}\left(\frac{1}{\prod_{i=1}^n c_i} - \sum_{i=1}^n\frac{q^{c_i}}{c_i\prod_{1\leq j\leq n, j\neq i} (c_j-c_i)}\right)$, assuming the entries of $C$ are all pairwise distinct.
\end{enumerate}

\end{lemma}

\begin{proof}

The claim that $V(A,B,q)$ is continuous in $A$, $B$, and $q$ is obvious from the integral representation of the volume.

We now prove that the second case above holds and use it to deduce the first case.
We proceed by induction on $n$.
For $n=1$, the result is clear:

$$ \int_{x_1=q^{\frac{1}{a_1}}}^1 x_1^{b_1}dx_1 = \frac{1}{b_1+1}\left(1-q^{\frac{b_1+1}{a_1}}\right)= \frac{1}{a_1}\left( \frac{1}{c_1} - \frac{q^{c_1}}{c_1}\right)$$

Assuming the result for $n$, we may compute the $n+1$ dimensional volume as follows:

$$ \int_{x_{n+1}=q^{\frac{1}{a_{n+1}}}}^1 \frac{1}{\prod_{i=1}^n a_i}\left( \frac{1}{\prod_{i=1}^n c_i} - \sum_{i=1}^n \frac{(qx_{n+1}^{-a_{n+1}})^{c_i}}{c_i\prod_{1\leq j\leq n, j\neq i} (c_j-c_i)}\right)x_{n+1}^{b_{n+1}} dx_{n+1}$$

$$ = \frac{1}{\prod_{i=1}^{n+1} a_i}\left(\frac{1}{\prod_{i=1}^{n+1} c_i} - \sum_{i=1}^n \frac{q^{c_i}}{c_i\prod_{1\leq j\leq n, j\neq i} (c_j-c_i)} - \frac{q^{c_{n+1}}}{\prod_{i=1}^{n+1}c_i} + \sum_{i=1}^n \frac{q^{c_{n+1}}}{c_i\prod_{1\leq j \leq n+1, j\neq i} (c_j-c_i)} \right)$$

By a classical theorem of Sylvester \cite{Sylvester}, we have that $1=\sum_{i=1}^{n+1}\prod_{j\neq i}^{n+1} \frac{c_j}{c_j-c_i}$, or equivalently $\frac{1}{\prod_{i=1}^{n+1}c_i}-\sum_{i=1}^n\frac{1}{c_i\prod_{1\leq j \leq n+1, j\neq i} (c_j-c_i)} = \frac{1}{c_{n+1}\prod_{i=1}^n(c_i-c_{n+1})}$, and so the coefficient of $q^{c_{n+1}}$ is as claimed.

We now show the first case.
By another classical theorem of Sylvester \cite{Sylvester}, we have that $h_{d-n+1}(x_1,\cdots,x_n)=\sum_{r=1}^n x_r^d \prod_{j\neq r}\frac{1}{x_r-x_j}$.
Assuming $c_i\neq c_j$ for $i\neq j$, let $0<q<1$, expand the $m$-dimensional volume function as a series in $\log(q)$, and apply the result of Sylvester to obtain that the volume function is equal to 

$$\frac{1}{\prod_{i=1}^n a_i}\sum_{i=0}^\infty \frac{h_{i-n}(c_1,\cdots,c_n)(\log(q))^i}{i!} $$

when $c_i\neq c_j$ for $i\neq j$.

By the root test and the estimate $|h_{i-n}(c_1,\cdots,c_n)|\leq \binom{i-1}{n-1}\max_j |c_j|^{i-n}$, this series converges for all $C$ (in particular, for all $A>0$ and all $B \geq 0$) and $0<q<1$ (even when $c_i=c_j$ for some $i\neq j$). 
Since the above series agrees with $\frac{1}{\prod_{i=1}^n a_i}\left( \frac{1}{\prod_{i=1}^n c_i} - \sum_{i=1}^n \frac{q^{c_i}}{c_i\prod_{1\leq j\leq n, j\neq i} (c_j-c_i)}\right)$ on the open dense subset where both are defined and $V(A,B,q)$ is continuous, we see that $V(A,B,q)$ is in fact analytic in $A$, $B$, and $\log(q)$ with a series representation given by $\frac{1}{\prod_{i=1}^n a_i}\sum_{i=0}^\infty \frac{h_{i-n}(c_1,\cdots,c_n)(\log(q))^i}{i!} $ for $A>0$, $B\geq 0$, and $0<q<1$.

\end{proof}

Since all partial derivatives of all orders of an analytic function defined on an open set $U$ are again analytic, convergent, and continuous on $U$, we see that the negative of the derivative of $V(A,B,q)$ with respect to $q$ is again analytic for all $0<q<1$ and all $a_i\neq 0$, and is given by the formula $$\frac{1}{\prod_{i=1}^n a_i}\sum_{i=0}^\infty \frac{h_{i-n+1}(c_1,\cdots,c_n)(\log q)^i}{i!}.$$

To show this function has a continuous extension to $q=0$ which is zero if all entries of $A$ are even integers while the conditions on $A$ and $B$ from Theorem \ref{MainTheoremLocal} are met, we analyze the function $$F(C,q)=\sum_{i=0}^\infty \frac{h_{i-n+1}(c_1-1,\cdots,c_n-1)(\log q)^i}{i!}.$$ 

\begin{lemma}
\label{FinalComputation}

Suppose $C=(c_1,\cdots,c_n)$ is a sequence of real numbers with either $c_i>1$ for all $i$ or $c_i\geq 1$ for all $i$ with exactly one index $i_0$ such that $c_{i_0}=1$.
Then $F(C,q)$ has continuous extension to $q=0$.
Additionally, the continuous extension takes the value $0$ if no $c_i=1$.

\end{lemma}

\begin{proof}

Up to relabeling, we may assume $c_1\leq c_2\leq \cdots \leq c_n$.
In the case where all entries of $C$ are pairwise distinct, we may apply formula (2) from Lemma \ref{VolumeComputation} to see that $$F(C,q)=\sum_{i=0}^n \frac{q^{c_i-1}}{\prod_{j\neq i} (c_j-c_i)},$$ which has limit $0$ as $q\to 0$ if $c_1>1$ and limit $\frac{1}{\prod_{j=2}^n (c_j-1)}$ if $c_1=1$. 

In the remaining case where some entries of $C$ are equal to each other, assume again that we have ordered $c_1\leq c_2\leq\cdots\leq c_n$.
Since $F(C,q)$ is continuous in $C$ for all $C$ and $0<q<1$, for any $0<q<1$ we may write $F(C,q)=\lim_{C'\to C} F(C',q)$ for $C'$ a sequence of real numbers such that all entries of $C'$ are distinct.

Let $Y_1=\{k\in\ZZ_{>0} \mid c_1=c_k\}$, and let $Y_l$ be inductively defined as $\{k\in \NN \mid c_{\min \ZZ_{>0}\setminus (Y_1\cup\cdots\cup Y_{l-1})} = c_k\}$.
Using the equality $F(C',q)=\sum_{i=0}^n \frac{q^{c_i-1}}{\prod_{j\neq i} (c_j-c_i)}$ and the linearity of the limit, we see that we can write $$\lim_{C'\to C} F(C',q) = \sum_{l} \lim_{C'\to C} \sum_{i\in Y_l} \frac{q^{c_i-1}}{\prod_{j\neq i} (c_j-c_i)}$$ assuming each limit $\lim_{C'\to C} \sum_{i\in Y_l} \frac{q^{c_i-1}}{\prod_{j\neq i} (c_j-c_i)}$ exists.
Factoring out $(c_j-c_i)$ for each $j\notin Y_l$, we see that it is enough to treat the case of $$\lim_{(c_1,\cdots,c_m)\to(c_1,c_1,\cdots,c_1)}  \sum_{i=0}^m \frac{q^{c_i-1}}{\prod_{j\neq i} (c_j-c_i)}.$$

Rewriting, we have $\lim_{(c_1,\cdots,c_m)\to(c_1,c_1,\cdots,c_1)} \sum_{i=0}^\infty \frac{h_{i-m+1}(c_1-1,\cdots,c_m-1)(\log q)^i}{i!}$.
Applying continuity of this power series in the $c_i$, we see that this is equal to $\sum_{i=0}^\infty \frac{h_{i-m+1}(c_1-1,\cdots,c_1-1)(\log q)^i}{i!}$.
As $h_{i-m+1}(c_1-1,\cdots,c_1-1)=\binom{i-m+1+m-1}{m-1}(c_1-1)^{i-m+1}=\binom{i}{m-1}(c_1-1)^{i-m+1}$, we may write this series as $$\sum_{i=m-1} \frac{(c_1-1)^{i-m+1}(\log q)^i}{(m-1)!(i-m+1)!}$$ and reindexing with $j=i-m+1$ gives $$\frac{(\log q)^{m-1}}{(m-1)!}\sum_{j=0}^\infty \frac{((c_1-1)\log q)^j}{j!} = q^{c_1-1}\frac{(\log q)^{m-1}}{(m-1)!}.$$

Since in this case $c_1>1$ by our assumptions on the permissible $C$, this has limit $0$ as $q\to 0$, and we have proven the lemma.

\end{proof}

If all nonzero entries of $A$ are even integers, then by the conditions in Theorem \ref{MainTheoremLocal}, all $c_i$ must be greater than $1$.
By Lemma \ref{ReductionToSmallBox}, we have shown that the conditions of Lemma \ref{ReductionToBox} are satisfied, and thus we have proven Theorem \ref{MainTheoremLocal}, which in turn implies Theorem \ref{MainTheorem} by the reductions in Sections 4 and 5.

\appendix

\section{Recollections in Complex-Analytic Geometry}

In these appendices, we state a collection of definitions and results which we need to reference in the main body of the text.
We abandon our convention that when talking about geometric spaces, Roman letters are reserved for schemes and script letters are reserved for complex-analytic spaces.
Each time we define a geometric space in the statement of a theorem or result, we will make it clear what we are referring to.

\subsection{Basic Definitions}

We deal with analytic spaces over $k=\RR$ or $\CC$ in this paper.

\begin{definition}
\label{AnalyticModel}

Let $\cO_{k^n}$ be the sheaf of analytic functions on $k^n$.
Let $U$ be an open connected subset of $k^n$, and fix finitely many analytic functions $f_1,\cdots,f_m\in\cO_{k^n}(U)$.
Let $V\subset U$ be the common vanishing locus of the functions $f_1,\cdots,f_m$.
Define a sheaf of rings $\cO_{V}$ on $V$ as the restriction of $\cO_U/(f_1,\cdots,f_m)$, where $\cO_U$ is the restriction to $U$ of the sheaf $\cO_{k^n}$.
We call the locally ringed space $(V,\cO_V)$ a model.

\end{definition}

\begin{definition}
\label{AnalyticSpace}

A $k$-analytic space $X$ is a locally ringed space $(X,\cO_X)$ so that every point $x\in X$ has an open neighborhood $U\subset X$ isomorphic as a locally ringed space to a model $V$ as described in Definition \ref{AnalyticModel}. 


\end{definition}

\begin{definition}
\label{AnalyticSpaceMorphism}

A morphism of $k$-analytic spaces is a morphism of locally-ringed spaces such that the induced map on structure sheaves is $k$-linear.

\end{definition}

\begin{definition}
\label{ComplexConjugation}

Suppose $(X,\cO_X)$ is a complex-analytic space.
By a complex conjugation, we mean a map of locally-ringed spaces $\sigma:(X,\cO_X)\to(X,\cO_X)$ which is an order-two involution on the topological space $X$ and $\sigma_U:\cO_X(\sigma(U))\to\cO_X(U)$ satisfies $z\sigma_U(s)=\sigma_U(\ol{z}s)$ for $s\in\cO_X(\sigma(U))$ and $z\in\CC$.

\end{definition}

\begin{lemma}
\label{ComplexConjugationLemma}

Suppose $X$ is a complex-analytic variety and $\sigma$ is a complex conjugation on $X$.
The following claims hold:

\begin{enumerate}
\item $(X^\sigma,\cO_X^\sigma)$ is a real-analytic variety.
\item If $x\in X^\sigma\subset X$ is a smooth point as a complex-analytic variety, then $x\in X^\sigma$ is a smooth point as a real-analytic variety.
\item If $x\in X^\sigma$ is a smooth point, the the real dimension of $X^\sigma$ at $x$ coincides with the complex dimension of $X$ at $x$.
\item For every point $x\in X^\sigma$, we may choose a $\sigma$-invariant open set $U_x$ so that $U_x$ is isomorphic to a local model inside $\CC^n$, where the ideal sheaf is generated by real functions.
Equivalently, we may choose a local model so that the isomorphism between $U_x$ and the local model is an intertwiner for the complex conjugation on $X$ and the complex conjugation on the local model inherited from $\CC^n$.
\end{enumerate}

\end{lemma}

\begin{proof}

See \cite{GMTRealAnalytic}, Chapter II Section 4.

\end{proof}

\subsection{Measures on Analytic Spaces}

In a manner similar to Definition \ref{AlgebraicMeasuresGeneralization}, we define how to construct a measure on the complex-conjugation invariant locus of a complex-analytic variety:

\begin{definition}
\label{ComplexConjugateMeasure}

Suppose $X$ is a Gorenstein complex-analytic space with a complex conjugation $\sigma$.
Let $\W_X$ be the dualizing sheaf of $X$.
For $\w\in\W_X(X)^\sigma$ a global section invariant under $\sigma$, we define a measure $|\w|$ on $X^\sigma$ as follows.
Recall that $\W_X|_{X^{sm}}$ is isomorphic to the line bundle of top differential forms on $X^{sm}$.
Given a relatively compact open subset $U\subset X^\sigma$ and an analytic diffeomorphism $\Psi$ between $U\cap X^{sm}$ and $W\subset \RR^n$, we may write $$\Psi^*\w = g dx_1\wedge \cdots\wedge dx_n$$ for some $g:W\to \RR$, and define $$|\w|(U) = \int_W |g|d\l$$ where $|g|$ is the usual absolute value ($g$ is a real-valued function by Lemma \ref{ComplexConjugationLemma}) and $\l$ is the standard Lebesgue measure on $\RR^n$.
By the change of variables formula, this definition is independent of the diffeomorphismn $\Psi$.
There is a unique extension of $|\w|$ to a (possibly infinite) Borel measure on $X^\sigma$, which we also denote $|\w|$. 

\end{definition}

\subsection{Stein Spaces}

A particularly useful concept in the study of complex-analytic spaces is that of Stein spaces.

\begin{definition}
\label{SteinDefinition}

Let $\cX$ be a second countable complex-analytic space, and $\cO_X$ its structure sheaf.
$X$ is called Stein if it satisfies the following criteria:

\begin{enumerate}
\item For every pair of distinct points $x,y\in X$, there exists a holomorphic function $f\in\cO_X(X)$ such that $f(x)\neq f(y)$.
\item Every local ring $\cO_{X,x}$ is generated by functions in $\cO_X(X)$.
\item $X$ is holomorphically convex - that is, for every compact $K\subset X$, the set $$ \widehat{K}_X = \{ p\in X : |f(p)|\leq \max_{x\in K} |f(x)| \forall f\in\cO_X(X)\}$$ is again compact.
\end{enumerate}

\end{definition}

\begin{theorem}(Cartan's Theorem B)
\label{CartanB}

Let $X$ be a complex-analytic space.
$X$ is Stein if and only if for any coherent analytic sheaf $\cF$ on $X$, $H^i(X,\cF)=0$ for all $i>0$.

\end{theorem}

\begin{theorem}(Cartan's Theorem A)
\label{CartanA}

Let $X$ be a Stein complex-analytic space.
Then every coherent sheaf $\cF$ on $X$ is generated by global sections.

\end{theorem}

\begin{proof}

See \cite{CartanTheoremB}.
Cartan's Theorem B is presented before Cartan's Theorem A because Theorem A is actually a consequence of Theorem B.

\end{proof}

\begin{lemma}
\label{SteinBasics}

Let $X,Y$ be complex-analytic spaces.

\begin{enumerate}
\item If $X\subset Y$ is a closed complex-analytic subvariety and $Y$ is Stein, then $X$ is also Stein.
\item If $X,Y$ are both Stein, then $X\times Y$ is also Stein.
\item If $X\subset \CC$ is open, $X$ is Stein.
\end{enumerate}

\end{lemma}

\begin{proof}

In each case, verifying that global holomorphic functions separate points and generate local rings at points is immediate.
We concentrate on proving holomorphic convexity.

(1): Assume $Y$ is holomorphically convex.
Let $K\subset X$ be a compact set and $x\in\widehat{K}_X$.
This means that $|f(x)|\leq \max_{p\in K} |f(p)|$ for all $f\in\cO_X(X)$.
In particular, this statement holds for every restriction of a function $f\in\cO_Y(Y)$ to $X$, so $x$ must also belong to $\widehat{K}_Y$.
Therefore $\widehat{K}_X\subset \widehat{K}_Y$.
Since $\widehat{K}_Y$ was assumed to be compact and $\widehat{K}_X$ is closed, $\widehat{K}_X$ is again compact and therefore $X$ is holomorphically convex.

(2): Assume $X,Y$ are both holomorphically convex.
Let $K\subset X\times Y$ be a compact set, and let $(x,y)\in \widehat{K}_{X\times Y}$.
Then $|f(x,y)|\leq \max_{p\in K} |f(p)|$ for all $f\in\cO_{X\times Y}(X\times Y)$, and in particular this holds for all functions which are constant in the $X$ or $Y$ direction.
If $f$ is constant in the $Y$ direction, this says that $|f(x)|\leq \max_{p\in pr_X K} |f(p)|$ for all $f\in\cO_X(X)$, so $x$ must belong to $\widehat{pr_X K}$.
A similar argument for functions which are constant in the $X$ direction shows that $y\in\widehat{pr_Y K}$.
Therefore $\widehat{K}_{X\times Y}\subset \widehat{pr_X K}\times \widehat{pr_Y K}$, so $\widehat{K}_{X\times Y}$ must again be compact as a closed subset of a compact set.

(3): We show directly that every open subset $X\subset \CC$ is holomorphically convex.
Let $K$ be a compact subset of $X$.
$\widehat{K}_{X}\subset X$ is closed in the induced topology, so to show that it is compact, it is enough to show that $\widehat{K}_X$ avoids the boundary of $X$ and is bounded.
To see that $\widehat{K}_X$ is bounded, consider the definition of holomorphic convexity applied to the function $z:\CC\to\CC$.
If $x\in \widehat{K}_X$, then $|x|\leq \max_{p\in K} |p|$, so $\widehat{K}_X$ is bounded.
To see that $\widehat{K}_X$ avoids the boundary, let $a\in\CC$ belong to $\pd X$.
Apply the definition of holomorphic convexity to $z\mapsto\frac{1}{z-a}$.
Since $K$ is compact and does not contain $a$, $|\frac{1}{z-a}|$ is bounded on $K$.
This means that $\widehat{K}_X$ cannot approach $a$, and we are done.

\end{proof}

\begin{corollary}
\label{SteinConsequences}

The following statements follow easily from Lemma \ref{SteinBasics}:

\begin{enumerate}
\item $\CC^n$ is Stein.
\item Any polydisc is Stein.
\item Any complex-analytic variety which admits a closed analytic embedding into $\CC^n$ is Stein.
\item If $X$ is an affine algebraic variety defined over a subfield $k\subset \CC$, then $X(\CC)$ is Stein.
\item Let $f:X\to Z$ and $g:Y\to Z$ be holomorphic and assume that $X,Y$ are Stein.
Then $X\times_Z Y$ is Stein.
\end{enumerate}

\end{corollary}

\begin{proof}

(1),(2): Apply statements (2) and (3) from Lemma \ref{SteinBasics}.

(3): Apply (1) from Lemma \ref{SteinBasics}.

(4): Since $X$ is affine, it admits a closed embedding $X\hookrightarrow \AA^N_k$ for some $N$.
Taking $\CC$ points, we see that $X(\CC)\hookrightarrow\CC^N$ and we are done by (3).

(5): $X\times_Z Y$ is a closed subvariety of $X\times Y$, so by applying statements (1) and (2) from Lemma \ref{SteinBasics}, we have the result.

\end{proof}

\subsection{Analytic Duality}

For some constructions related to rational singularities on complex-analytic varieties, it is important that we develop a complex-analytic analogue of portions of Grothendieck duality.
In this subsection, we will define the dualizing complex and recall a duality theorem of Ramis, Ruget, and Verdier.

\begin{definition}
\label{DualizingComplexAnalytic}

Let $X$ be a complex-analytic space of dimension $n$.
Let $U$ be an open subset of $X$ isomorphic to a model in an open subset $Z\subset \CC^m$.
We define $$\W_X^\bullet|_U = R\Hom_X(\cO_X,\W_Z[m])|_U$$ where $\W_Z$ is the line bundle of top differential forms on $Z$, which is the dualizing complex on $Z$ up to a shift of $m$.

\end{definition}

\begin{lemma}
\label{DualityLemmaAnalytic}

Let $X$ be a complex-analytic space.
Define an endofunctor $D_X$ on the category of complexes of sheaves on $X$ with coherent cohomology by $D_X(\cF^\bullet)=R\Hom_X(\cF^\bullet,\W_X^\bullet)$.
Then the canonical morphism from $\cF^\bullet$ to $D_XD_X\cF^\bullet$ is a quasi-isomorphism.

\end{lemma}

\begin{proof}

See \cite{RamisRuget}.

\end{proof}

\begin{theorem}
\label{DualityTheoremAnalytic}

Suppose $f:X\to Y$ is a proper morphism of analytic spaces.
Let $\cF^\bullet$ be a complex of sheaves on $X$ with coherent cohomology, and let $\cG^\bullet$ be a bounded-below complex of sheaves on $Y$ with coherent cohomology.
There is a canonical isomorphism $$Rf_*R\Hom_X(\cF^\bullet,R\Hom_X(Lf^*R\Hom_Y(\cG^\bullet,\W_Y^\bullet),\W_X^\bullet)) \to R\Hom_Y(Rf_*\cF^\bullet,\cG^\bullet).$$

\end{theorem}

\begin{proof}

See \cite{RamisRugetVerdier}.

\end{proof}

\begin{corollary}
\label{EasierDualityAnalytic}

Suppose $f:X\to Y$ is a proper morphism of analytic spaces.
Let $\cF^\bullet$ be a complex of sheaves on $X$ with coherent cohomology.
There is a canonical isomorphism $$Rf_*R\Hom_X(\cF^\bullet,\W_X^\bullet)\to R\Hom_Y(Rf_*\cF^\bullet,\W_Y^\bullet).$$

\end{corollary}

\begin{proof}

Applying Theorem \ref{DualityTheoremAnalytic} with $\cG^\bullet=\W_Y^\bullet$, we see that we have a canonical isomorphism $$Rf_*R\Hom_X(\cF^\bullet,R\Hom_X(Lf^*R\Hom_Y(\W_Y^\bullet,\W_Y^\bullet),\W_X^\bullet)) \to R\Hom_Y(Rf_*\cF^\bullet,\W_Y^\bullet).$$
Noting that $R\Hom_Y(\W_Y^\bullet,\W_Y^\bullet)=\cO_Y$, $Lf^*\cO_Y=\cO_X$, and that $R\Hom_X(\cO_X,-)$ is the identity functor, we obtain the result.

\end{proof}

\section{Resolution of Singularities}

\subsection{Algebraic Resolution of Singularities}

We will use Hironaka's theorem on resolution of singularities in characteristic zero, found in \cite{HironakaResolutionReference}.
The algebraic machinery introduced in this section is primarily intended as background, and may be helpful in understanding resolution of singularities in the analytic context in the next subsection.

\begin{definition}
\label{AlgebraicResolution}

Let $X$ be an algebraic variety.

\begin{itemize}
\item A resolution of singularities of $X$ is a proper map $\pi: Y\to X$ such that $Y$ is smooth and $\pi$ is a birational equivalence.
\item A strong resolution of singularities of $X$ is a resolution which is an isomorphism over the smooth locus of $X$.
\item A subvariety $D\subset X$ is said to be a normal crossings divisor if for any $x\in D$ there exists an etale neighborhood $\phi:U\to X$ of $x$ and an etale map $\alpha:U\to \AA^n$ such that $\phi^{-1}(D)=\a^{-1}(D')$, where $D'\subset \AA^n$ is a union of coordinate hyperplanes.
\item A subvariety $D\subset X$ is said to be a strict normal crossings divisor if for any $x\in D$ there exists a Zariski neighborhood $U\subset X$ of $x$ and an etale map $\a:U\to\AA^n$ such that $D\cap U=\a^{-1}(D')$, where $D'\subset \AA^n$ is a union of coordinate hyperplanes.
\item We say that a resolution of singularities $\pi:\widetilde{X}\to X$ resolves (respectively strictly resolves) a closed subvariety $D\subset X$ if $\pi^{-1}(D)$ is a normal crossings divisor (respectively a simple normal crossings divisor).
In this case we will also say that $\pi:\widetilde{X}\to X$ is a resolution (respectively a strict resolution) of the pair $(X,D)$.
\item Let $D\subset X$ be a subvariety of codimension 1. 
A strong resolution of the pair $(X,D)$ is a strict resolution $\pi:\widetilde{X}\to X$ which is an isomorphism outside of the singular locus of $X$ and the singular locus of $D$.
\item Let $\pi:\widetilde{X}\to X$ be a resolution of singularities.
Let $U\subset X$ be the maximal open set on which $\pi$ is an isomorphism.
Let $D\subset X$ be a subvariety. The strict transform of $D$ is defined to be $\ol{\pi^{-1}(D\cup U)}$.
\end{itemize}

\end{definition}

\begin{theorem}[Hironaka]
\label{HironakaAlgebraicResolution}

Let $D\subset X$ be a pair of algebraic varieties.
Assume that $D\subset X$ is of codimension 1 and let $U\subset X$ be a smooth open subset such that $U\cap D$ is also smooth.
Then there exists a resolution of singularities $\pi:\widetilde{X}\to X$ that resolves $D$ and such that $\pi:\pi^{-1}(U)\to U$ is an isomorphism.

\end{theorem}

There is a standard method to resolve a normal crossings divisor to a strict normal crossings divisor.
This implies the following corollary.

\begin{corollary}

Any pair of algebraic varieties $D\subset X$ admits a strong resolution.

\end{corollary}

The following lemma is standard.

\begin{lemma}
\label{EmbeddedPairAlgebraicResolution}

Let $D\subset X$ be a pair of algebraic varieties such that $D$ is irreducible, has codimension one in $X$, and is not contained in the singular locus of $X$.
Let $\pi:\widetilde{X}\to X$ be a strong resolution of the pair $(X,D)$ and let $D'$ be the strict transform of $D$.
The $\pi|_{D'}:D'\to D$ is a resolution of singularities.

\end{lemma}

\subsection{Analytic Resolution of Singularities}

We use Hironaka's theorem for resolution of singularities for analytic spaces, found in \cite{HironakaResolutionReference}.
A more modern overview may be found in \cite{Wlodarczyk}.

\begin{definition}
\label{AnalyticResolution}

Let $X$ be a complex-analytic variety.

\begin{itemize}
\item An analytic subvariety $D\subset X$ is said to be a strict normal crossings divisor if for any $x\in D$ there exists a neighborhood $U\subset X$ and an open set $V\subset \CC^n$ such that there is an isomorphism $\a:U\cong V$ and $D=\a^{-1}(D')$, where $D'\subset V$ is a union of coordinate hyperplanes.
\item A resolution of singularities of $X$ is a manifold $\widetilde{X}$ with a proper bimeromorphic map $\pi:\widetilde{X}\to X$ such that $\pi$ is an isomorphism over the nonsingular part of $X$, the inverse image of the singular locus is a strict normal crossings divisor, and $\pi$ is functorial with respect to local analytic isomorphisms.
\item Let $\pi:\widetilde{X}\to X$ be a resolution of singularities.
Let $U\subset X$ be the maximal open set on which $\pi$ is an isomorphism.
Let $D\subset X$ be a subvariety. The strict transform of $D$ is defined to be $\ol{\pi^{-1}(D\cup U)}$.
\item Let $X$ be a smooth complex-analytic space (i.e. a manifold) and $D\subset X$ be an analytic subspace.
There exists a manifold $\widetilde{X}$ together with a proper bimeromorphic map $\pi:\widetilde{X}\to X$ and strict normal crossings locally finite divisor $E\subset X$ such that the strict transofrm $\widetilde{D}$ of $D$ is smooth and has simple normal crossings with $E$.
Furthermore, the support of $E$ is the exceptional locus of $\pi$, and $\pi$ factors as a sequence of blowups in smooth centers.
\item Let $D\subset X$ be a subvariety of codimension 1. 
A strong resolution of the pair $(X,D)$ is a resolution $\pi:\widetilde{X}\to X$ which is an isomorphism outside of the singular locus of $X$ and the singular locus of $D$.
\end{itemize}

\end{definition}

\begin{theorem}[Hironaka]
\label{HironakaAnalyticResolution}

Each sort of resolution defined above in definition \ref{AnalyticResolution} exists.

\end{theorem}

An analogue of \ref{EmbeddedPairAlgebraicResolution} holds in the analytic situation.

\subsection{Rational Singularities}

Here we recall the definition of rational singularities for both algebraic varieties and complex-analytic varieties.
Both definitions depend on the notion of resolution of singularities, defined previously in the appendix.

\begin{definition}
\label{RationalSingularitiesDefinition}

We define rational singularities for both algebraic and analytic spaces.

\begin{enumerate}
\item We say that an algebraic variety or complex-analytic variety $X$ has rational singularities if for any (equivalently, some) resolution of singularities $\pi: \widetilde{X}\to X$, the natural map $\cO_X\to R\pi_*(\cO_{\widetilde{X}})$ is an isomorphism.
\item Suppose $X$ is an algebraic variety.
A point $x\in X(k)$ is said to be a rational singularity if there is a Zariski neighborhood of $x$ that has rational singularities.
\item Suppose $X$ is a complex-analytic variety.
A point $x\in X$ is said to be a rational singlarity if there is a neighborhood of $x$ with rational singularities
\end{enumerate}

\end{definition}

\begin{theorem}
\label{MoreRationalSingularities}

Let $X$ be an algebraic variety.
The following are equivalent:

\begin{enumerate}
\item $X$ has rational singularities.
\item For any (equivalently, for some) resolution of singularities $\pi:\widetilde{X}\to X$, the trace map $Tr_\pi:R\pi_*\W_{\widetilde{X}}\to \W_X$ is an isomorphism.
\item $X$ is Cohen-Macaulay, and for any (equivalently, for some) resolution of singularities $\pi:\widetilde{X}\to X$, the trace map $Tr_\pi:R\pi_*\W_{\widetilde{X}}\to \W_X$ is an isomorphism.
\item $X$ is Cohen-Macaulay, and for any (equivalently, for some) resolution of singularities $\pi:\widetilde{X}\to X$, the trace map $Tr_\pi:R\pi_*\W_{\widetilde{X}}\to \W_X$ is onto.
\item $X$ is Cohen-Macaulay, normal, and for any (equivalently, for some) resolution of singularities $\pi:\widetilde{X}\to X$ the composition $R\pi_*\W_{\widetilde{X}}\stackrel{Tr_{\pi}}{\to} \W_X \to i_*(\W_{X^{sm}})$, where $i:X^{sm}\to X$ is the embedding of the smooth locus, is an isomorphism.
\end{enumerate}

If $X$ is affine, these conditions are also equivalent to the following:

\begin{enumerate}
\setcounter{enumi}{5}

\item $X$ is Cohen-Macaulay, and for any (equivalently, for some) strong resolution of singularities $\pi:\widetilde{X}\to X$ and any section $\w\in \W_X(X)$, there exists a top differential form $\widetilde{\w}\in \W_{\widetilde{X}}(\widetilde{X})$ agreeing with $\w$ on the smooth locus of $X$, i.e.

$$\w|_{X^{sm}} = \widetilde{\w}|_{X^{sm}}$$

where we consider $X^{sm}$ as a subset of both $X$ and $\widetilde{X}$.
\item $X$ is Cohen-Macaulay, normal, and for any (equivalently, for some) strong resolution of singularities $\pi:\widetilde{X}\to X$ and any top differential form $\w\in \W_X(X^{sm})$, there exists a top differential form $\widetilde{\w}\in \W_{\widetilde{X}}(\widetilde{X})$ agreeing with $\w$ on the smooth locus of $X$, i.e.

$$\w = \widetilde{\w}|_{X^{sm}}$$

where we consider $X^{sm}$ as a subset of both $X$ and $\widetilde{X}$.
\end{enumerate}

\end{theorem}

\begin{proof}

See \cite{AvniAizenbud}, Proposition B.7.2.

\end{proof}

\begin{theorem}
\label{MoreRationalSingularitiesAnalytic}

Suppose $X$ is a complex-analytic variety.
The following are equivalent:


\begin{enumerate}
\item $X$ has rational singularities.
\item For all (equivalently, for some) resolution of singularities $\pi:\widetilde{X}\to X$ the natural map $j:R^\bullet\pi_*\W_{\widetilde{X}}\to \W_X^\bullet$ is an isomorphism.
\item $X$ is Cohen-Macaulay, normal, and for all (equivalently, for some) resolution of singularities $\pi:\widetilde{X}\to X$ the natural map $j:\pi_*\W_{\widetilde{X}}\to \W_X$ is an isomorphism.
\item $X$ is Cohen-Macaulay, normal, and for all (equivalently, for some) resolution of singularities $\pi:\widetilde{X}\to X$ the natural map $j:\pi_*\W_{\widetilde{X}}\to \W_X$ is surjective.
\item Assume $X$ is also Stein.
$X$ is Cohen-Macaulay, normal, and for any section $\w\in\W_X$ and strong resolution of singularities $\pi:\widetilde{X}\to X$, there exists a top form $\widetilde{\w}\in\W_{\widetilde{X}}$ so that $\w|_{X^{sm}}=\widetilde{\w}|_{X^{sm}}$.
\end{enumerate}
\end{theorem}

\begin{proof}

(1)$\Leftrightarrow$(2): Let $j:R\pi_*\W_{\widetilde{X}}[\dim X]\to \W_X^\bullet$ be the natural map obtained by applying $D_X=R\Hom_X(-,\W_X^\bullet)$ to the map $\pi^{\#}:\cO_X\to R\pi_*\cO_{\widetilde{X}}$, where we rely on Theorem \ref{DualityTheoremAnalytic} for the identification of $D_X(R\pi_*\cO_{\widetilde{X}})\cong R\pi_*\W_{\widetilde{X}}$. 
Construct distinguished triangles $\cO_X\to R\pi_*\cO_{\widetilde{X}}\to M_X^\bullet \to \cO_X[1]$ and $R\pi_*\W_{\widetilde{X}}[\dim X]\to \W_X^\bullet \to N_X^\bullet \to R\pi_*\W_{\widetilde{X}}[1]$ in the derived category of complexes with coherent cohomology on $X$.
By Lemma \ref{DualityLemmaAnalytic}, $M_X^\bullet=D_X(N_X^\bullet)$ and $N_X^\bullet=D_X(M_X^\bullet)$.
Note that $R\pi_*\W_{\widetilde{X}}$ is concentrated in a single degree by Lemma \ref{AnalyticGR}, and $\cO_X$ is also concentrated in a single degree.
If $\pi^{\#}$ is an isomorphism, then $M_X^0=0$, which implies that $N_X^{-\dim X}=0$, and thus $j$ must also be an isomorphism.
Similarly, if $j$ is an isomorphism, then $\pi^{\#}$ must also be an isomorphism.

(2)$\Leftrightarrow$(3): This follows from Lemma \ref{AnalyticGR}, an analytic analogue of the Grauert-Riemenschneider Theorem.

(3)$\Leftrightarrow$(4): We need to prove that under the hypotheses that $X$ is Cohen-Macaulay and normal, the map $j:R^\bullet\pi_*\W_{\widetilde{X}}\to \W_X^\bullet$ is a monomorphism.
Let $U\subset X$ be an open set.
Because $X$ is Cohen-Macaulay and $\pi$ is a resolution of singularities, the map reduces to $j:\pi_*\W_{\widetilde{X}}\to \W_X$, a map of sheaves.
Suppose $\w\in \W_{\widetilde{X}}(\pi^{-1}(U)) = \pi_*\W_{\widetilde{X}}(U)$ is in the kernel.
Let $W\subset X$ be the open dense locus where $\pi$ is an isomorphism.
Then $\w|_{W\cap U}=0$, so $\w\in\W_{\widetilde{X}}(\pi^{-1}(U))$ is a holomorphic section which is zero on a set whose complement has codimension at least two.
By the Riemann Extension Theorem, $\w$ must be zero, and we have proven the claim.

(3)$\Rightarrow$(5): Since $j:\pi_*\W_{\widetilde{X}}\to\W_X$ is an isomorphism, we just compute the preimage of $\w$.
Since $j$ is the identity on $X^{sm}$, the claim follows.

(5)$\Rightarrow$(4): Since $X$ is Stein, the map $j$ is determined by its action on global sections by Cartan's Theorem A (Theorem \ref{CartanA}).
It suffices to show that $j(\widetilde{\w})=\w$.
By construction of $j$, $j(\widetilde{\w})|_{X^{sm}}=\w|_{X^{sm}}$, so $j(\widetilde{\w})-\w$ is zero on a set whose compliment has codimension at least two, and is thus zero.

\end{proof}

\begin{definition}
\label{FRSMore}

Suppose $X,Y$ are either algebraic varieties defined over a field $k$ of characteristic zero or complex-analytic varieties.
Suppose further that $Y$ is smooth.
Let $\varphi:X\to Y$ be a morphism.
We say that $\varphi$ is (FRS) if it is flat and:
\begin{enumerate}
\item in the case that $Y$ is an algebraic variety with $y\in Y(\ol{k})$, the fiber $X\times_Y y$ is reduced with rational singularities
\item in the case that $Y$ is a complex-analytic space with $y\in Y$, the fiber $X\times_Y y$ is reduced with rational singularities.
\end{enumerate}

\end{definition}

\end{document}